\documentclass[11pt]{article}
\usepackage[letterpaper,margin=1in]{geometry}
\usepackage[T1]{fontenc}
\usepackage{lmodern,amsmath,amssymb,amsthm,mathtools,microtype}
\usepackage[hidelinks]{hyperref}
\usepackage{fancyhdr}
\numberwithin{equation}{section}
\newtheorem{theorem}{Theorem}[section]
\newtheorem{proposition}[theorem]{Proposition}
\newtheorem{lemma}[theorem]{Lemma}
\newtheorem{corollary}[theorem]{Corollary}
\theoremstyle{remark}
\newtheorem{remark}[theorem]{Remark}
\newcommand{\R}{\mathbb R}
\newcommand{\Q}{\mathbb Q}

\newcommand{\Rea}{\operatorname{Re}}

\newcommand{\supp}{\operatorname{supp}}
\newcommand{\GL}{\operatorname{GL}}
\newcommand{\av}[1]{\left\langle #1\right\rangle_h}
\newcommand{\wt}{\widetilde}

\allowdisplaybreaks[2]
\title{Probing automorphic induction\\through zero correlations and secondary arithmetic terms}
\author{Yangbo Ye\\[4pt]
\small Department of Mathematics, The University of Iowa\\
\small Iowa City, Iowa 52242, USA\\
\small \texttt{yangbo-ye@uiowa.edu}}
\date{Expanded version --- October 22, 2026}
\hypersetup{pdftitle={Probing automorphic induction through zero correlations and secondary arithmetic terms},pdfauthor={Yangbo Ye},pdfsubject={Zero correlations and automorphic-induction probes}}

\begin{document}
\maketitle
\begin{abstract}
Zero correlations offer a way to test the cuspidal structure of an
automorphic $L$-function, but their leading self statistics retain
only coarse degree and multiplicity data. We therefore adjoin a
specified cuspidal $L$-function as a probe and study the resulting raw
mixed correlations. Assuming a factorization into cuspidal
$L$-functions over $\Q$, we show that a mixed pair statistic detects
the probe's occurrence and recovers its multiplicity. A suitable
Gevrey-order-two test eliminates the smooth-density terms and gives
a subexponential error. For arbitrary smooth plateau tests, the
prime number theorem of Harcos and Thorner gives an error smaller
than every fixed inverse power of $\log T$, without a
self-contragredience assumption. The stronger prime number theorem
of Liu and Ye gives a subexponential remainder under its stated
structural assumptions, with the possible norm-twist terms made
explicit. To obtain arithmetic information beyond this
overlap, we compute the correlation expansion at every fixed level
and inverse-logarithmic order within the original strict
small-support range. We treat the pair and higher correlations separately,
using Jiang's rational-field fixed-data estimates to justify the
required prime-power convergence. The pair calculation identifies
the mixed arithmetic finite parts and gives detailed error estimates,
including the contribution of norm twists. The higher-level
calculation retains connected prime-power invariants, whose first
cubic contribution appears at order $(\log T)^{-3}$. Explicit third-
and fourth-level formulas and examples explain the arithmetic
content of these coefficients and the limitations of self-correlations.
In particular, two products of Dirichlet $L$-functions have identical
gamma factors and all self-pair arithmetic coefficients but different
cubic constants. Varying powers of the probe separates
their mixed components. These formulas, together with the
distinct-index sieve, provide a collection of arithmetic probe
responses toward the further problem of recovering automorphic
induction structure.
\end{abstract}

\clearpage\begingroup\small\tableofcontents\endgroup\clearpage
\section{Introduction}

This paper has two related purposes. We compute secondary arithmetic
terms in the raw mixed $n$-level correlations of zeros of automorphic
$L$-functions, and use these correlations to probe the cuspidal
constituents of an automorphic induction. The second purpose gives a
reason to retain the secondary terms: a test $L$-function supplies
arithmetic information about a given $L$-function that its leading
self-correlation does not determine. We regard the probe developed
here as a starting point for further investigation.

Suppose that $\Pi$ is a unitary cuspidal representation of $\GL_m$ over
a number field $E$ of degree $\ell$ over $\Q$, and that
\begin{equation}\label{eq:induction}
 F(s)=L(s,\Pi)=\prod_{\alpha=1}^q
 L(s,\sigma_\alpha)^{r_\alpha},
 \qquad \sum_{\alpha=1}^q r_\alpha d_\alpha=\ell m.
\end{equation}
Here the $\sigma_\alpha$ are inequivalent unitary cuspidal
representations of $\GL_{d_\alpha}(\mathbb A_\Q)$, and
$r_\alpha\ge1$. Such a factorization is the structure we seek to
detect. Its existence is assumed in the application; we do not prove
automorphic induction over an arbitrary extension. No Galois or
prime-degree assumption on $E/\Q$ is used once
\eqref{eq:induction} is given. For the theory of base change and
automorphic induction in the established settings, see
Arthur--Clozel \cite{AC}.

Write $D_F=\ell m$ and $B_F=\sum_\alpha r_\alpha^2$.
The leading raw self-correlations in the small-support range contain
these two aggregate quantities. They do not by themselves give the
individual degrees $d_\alpha$ or multiplicities $r_\alpha$.
Now choose a unitary cuspidal representation $\pi$ of
$\GL_{m_0}(\mathbb A_\Q)$ and put $P(s)=L(s,\pi)$.
The zero multiset of $FP$ is the sum of those of $F$ and $P$.
Write $\mathcal R_{F,G}(T)$ for the raw mixed pair correlation defined
in \eqref{eq:correlation}. With the same normalization and test functions
for all four correlations, the exact difference
\[
 \mathcal R_{FP,FP}-\mathcal R_{F,F}-\mathcal R_{P,P}
       =2\mathcal R_{F,P}
\]
isolates the mixed pair correlation. We prove that one can choose a
Fourier test vanishing near zero for which
\begin{equation}\label{eq:introprobe}
 \frac{\mathcal R_{FP,FP}(T)-\mathcal R_{F,F}(T)
                   -\mathcal R_{P,P}(T)}{2J_\phi}
 =r_\Pi(\pi)+O\!\left(e^{-a\sqrt{\log T}}\right),
 \qquad J_\phi\ne0.
\end{equation}
Here $r_\Pi(\pi)$ is the multiplicity of $\pi$ in
\eqref{eq:induction}, with value zero if it is absent. The constant
$a>0$ may depend on the fixed representations and test. Thus a
positive response identifies a constituent of the specified degree
$m_0$ and gives its multiplicity. This test does not require subtraction
of a conductor or gamma-factor contribution. A second test, nonzero
near zero, determines the real part of a mixed arithmetic finite part
$C_{F,P}$. Section~\ref{sec:probes} gives these statements and their
proofs, as well as a finite-degree stopping criterion for a supplied
list of probes.

The study of local zero correlations has its classical starting point in
Montgomery's pair-correlation theorem and conjecture for the Riemann zeta
function \cite{Montgomery}.  His theorem, under the Riemann hypothesis and
with a restriction on the Fourier transform of the test function, connects
zero differences with the unitary random-matrix prediction.  Odlyzko's
computations \cite{Odlyzko} supplied extensive numerical evidence for this
connection and also exhibited arithmetic effects at finite height.
Hejhal \cite{Hejhal} treated triple correlations, and Rudnick and
Sarnak \cite{RS} developed the small-support theory of higher correlations
for principal $L$-functions.  The work of Katz and Sarnak
\cite{KSbook,KSsurvey} places these questions in a broader framework involving
Frobenius eigenvalues, monodromy and classical compact groups.  Their
framework also distinguishes local statistics high on the critical line
from the statistics of low-lying zeros in families.  The latter direction
includes the work of Iwaniec, Luo and Sarnak \cite{ILS}.  The present paper
concerns height averages of fixed $L$-functions.

Arithmetic terms beyond the universal limit also have an established
history.  Bogomolny and Keating \cite{BK95,BK96} investigated higher
correlations of the Riemann zeros by methods motivated by spectral
statistics.  Bogomolny, Bohigas, Leboeuf and Monastra \cite{BBLM} studied
finite-height corrections to the spacing distribution.  Conrey and Snaith
\cite{CS,CStriple,CSn} developed the ratios-conjecture approach and, assuming
that conjecture, obtained lower-order terms for triple and arbitrary-order
correlations of the Riemann zeros.  Rodgers \cite{Rodgers} proved a refined
pair-correlation formula for a class of smooth tests under the Riemann
hypothesis.  These results show why secondary arithmetic terms deserve to
be retained.  Our purpose is to compute them for raw mixed correlations of
fixed automorphic products in the original small-support range, and to use
the mixed coefficients to test for specified cuspidal constituents.

A related arithmetic question asks what information about an $L$-function
can be recovered from differences of its zeros or from correlations with
another zero set.  Ford and Zaharescu \cite{FZ} established such phenomena
for unnormalized zero differences, including an application to ranks of
elliptic curves.  Our scaling and test functions are different, but this
work is relevant to the use of a second $L$-function as a probe.
The superposition and functoriality problems studied by Liu and Ye
\cite{LY02,LY08,LY09} provide the direct starting point here.  We make the
probe explicit: adjoining $L(s,\pi)$ produces mixed correlations whose
leading overlap coefficient detects the multiplicity of the specified
representation $\pi$, while their secondary coefficients retain further
Rankin--Selberg and higher prime-power information.  This constituent test
is one step toward detecting automorphic induction; it is not a general
reconstruction theorem for an unknown induction decomposition.

The long-term goal is to use zero correlations to study the relations
among unramified Satake parameters that express automorphic induction.
Strong multiplicity one gives uniqueness from these local data, but
does not by itself provide a procedure for extracting the globally
compatible cuspidal constituents. Liu and Ye \cite{LY09} already
studied zero correlations for automorphic $L$-functions over $E$
with functoriality in view. The present probe develops this approach
by comparing a given zero set with those of auxiliary $L$-functions;
its secondary terms provide further arithmetic observables. Here
\eqref{eq:induction} is assumed. Identifying a sufficiently rich
family of probes and recovering the induction structure from their
responses are goals for subsequent work.

Our analytic calculation is carried out for arbitrary fixed finite
products over $\Q$, so the same formulas apply before a particular
automorphic-induction problem is chosen. At every fixed level and
every fixed inverse-logarithmic order, the raw mixed correlation has
an expansion whose coefficients are archimedean terms, Laurent
coefficients of mixed Rankin--Selberg coefficient series, and moments
of higher prime-power cumulants. The arithmetic at ramified primes
uses the actual standard local factors. Their difference from the
genuine ramified Rankin--Selberg factors affects the finite parts and
cannot be discarded.

The first connected contribution involving three coordinates occurs
at order $(\log T)^{-3}$. It comes from relations
$p^a p^b=p^{a+b}$, including $p\,p=p^2$, and is missed if one keeps
only pair matchings. We compute the third- and fourth-level formulas
explicitly. An example using four primitive Dirichlet characters
modulo $13$ gives two products with identical gamma factors and all
self-pair arithmetic coefficients, but different cubic constants.
This proves a strict gain in arithmetic information at level three.
Applying the same expansion to $FP^q$ and varying the nonnegative
integer $q$ separates polarized higher invariants by polynomial
interpolation. These are further observables for the probe, although
they do not yet give a general reconstruction theorem.

The arithmetic inputs used here are established in the stated scope.
The rational-field fixed-data estimates in Jiang's preprint
\cite{Jiang} imply both Hypothesis H and the stronger mixed
power-saving tail needed below. Section~\ref{sec:established}
specifies the input and a harmless logarithmic correction in its use.
Harcos--Thorner \cite[Theorem 2.1]{HT}, applied with modulus $1$,
supplies a prime number theorem with an arbitrarily large fixed
logarithmic saving for every cuspidal cross pair. Thus every smooth
plateau gives a pair-correlation remainder smaller than each fixed
inverse power of $\log T$, without self-contragredience. Theorem~2.3
of Liu--Ye \cite{LY07}, with the corrected local treatment in
\cite{Ycorr}, supplies the stronger subexponential prime number
theorem when each cross pair contains a self-contragredient member.
Together with Jiang's prime-power tail, these prime number theorems
also imply the corresponding weighted prime Mertens formulas and
Selberg orthogonality for pairs unrelated by a norm twist.
For Gevrey-order-two Fourier tests constant near zero, we prove the
subexponential pair refinement for arbitrary fixed products from
local meromorphy. The general higher-level expansion already applies
to arbitrary fixed products and uses local continuation rather than
either quantitative prime number theorem. These results keep the
original strict Fourier-support range.

The distinction between raw and distinct-index correlations is
maintained throughout. We give the exact partition sieve, including
the changes in the cutoff and test under a contraction, and show
that it preserves the secondary expansion. The supplied probe detects
a candidate constituent; finding a sufficiently rich collection of
probes remains a further problem. The complete self-correlation
hierarchy is not an injective invariant in the generality considered
here, as Section~\ref{sec:information} explains.

The expanded version retains the separate pair calculation, the
PNT-based treatment of norm twists, the fourth-level formula, and the
examples. The shorter journal version combines the pair and higher
calculations into a single argument and omits the examples. All
representations, correlation levels, truncation orders, and test
functions are fixed as $T\to\infty$. Constants may depend on these
data and on nonzero norm-twist separations; no uniformity in a twist
tending to zero, or interchange with a limit $n\to\infty$, is asserted.

\section{Notation and main results}

Let $\mathcal S$ be a finite set of inequivalent unitary cuspidal
representations over $\Q$. Write
\begin{equation}\label{eq:products}
 F(s)=\prod_{\sigma\in\mathcal S}L(s,\sigma)^{r_\sigma},
 \qquad G(s)=\prod_{\sigma\in\mathcal S}L(s,\sigma)^{s_\sigma},
\end{equation}
where the exponents are nonnegative integers and neither product is
empty. Put
\begin{equation}\label{eq:degrees}
 D_F=\sum_\sigma r_\sigma d_\sigma,\quad
 D_G=\sum_\sigma s_\sigma d_\sigma,\quad
 r_{F,G}=\sum_\sigma r_\sigma s_\sigma.
\end{equation}
In particular, $r_{F,F}=B_F:=\sum_\sigma r_\sigma^2$.

At every finite prime use the actual standard local factors, with zero
roots appended when necessary. Thus
\begin{equation}\label{eq:coefficients}
 -\frac{F'}F(s)=\sum_{p,k\geq1}
  \frac{b_F(p^k)\log p}{p^{ks}},\qquad
 b_F(p^k)=\sum_\sigma r_\sigma\sum_{j=1}^{d_\sigma}
                    \alpha_\sigma(p,j)^k
\end{equation}
in a right half-plane. Define $b_G$ similarly. When a sum is written over
integers, $b_F(n)$ and $b_G(n)$ may be set to zero off prime powers.
Write
\[
 w_{F,G}(p^k)=b_F(p^k)\overline{b_G(p^k)}.
\]

The completed function is
\begin{equation}\label{eq:completion}
 \Lambda_F(s)=Q_F^{s/2}\prod_{j=1}^{D_F}
       \Gamma_\R(s+\mu_{F,j})F(s),\qquad
 \Gamma_\R(s)=\pi^{-s/2}\Gamma(s/2).
\end{equation}
The zero multiset used below consists of zeros of this completed
function, counted with multiplicity, and does not include its poles.
In particular, zeros of the finite $L$-function canceled by gamma poles
are not counted. Set $\gamma_F=(\rho_F-1/2)/i$ for such a zero $\rho_F$.
The ordinates may be complex. No Riemann hypothesis is assumed.

Fix $c>0$ and set
\begin{equation}\label{eq:scales}
 L=\log T,\qquad U=T/c,\qquad \lambda=cL,
 \qquad U\lambda=TL.
\end{equation}
Choose a nonzero real even function $g\in C_c^\infty(\R)$, with
$\supp g\subset[-A,A]$, and put
\begin{equation}\label{eq:h}
 h(t)=\int_\R g(v)e^{itv}\,dv,\qquad
 \kappa=\int_\R h(t)^2\,dt=2\pi\int_\R g(v)^2\,dv>0.
\end{equation}
Let $\phi\in C_c^\infty(\R)$ be real and even, and assume
\begin{equation}\label{eq:support}
 \supp\phi\subset[-b,b],\qquad cb<1.
\end{equation}
Define
\[
 \psi(x)=\int_\R\phi(u)e^{2\pi iux}\,du,\qquad
 J_\phi=\int_\R |u|\phi(u)\,du.
\]
The normalized raw mixed pair correlation is
\begin{equation}\label{eq:correlation}
 \mathcal R_{F,G}(T)=\frac{2\pi}{\kappa TL}
 \sum_{\gamma_F,\gamma_G}
 h(\gamma_F/U)h(\gamma_G/U)
 \psi\!\left(\frac{\lambda(\gamma_F-\gamma_G)}{2\pi}\right).
\end{equation}
All pairs are retained, including coincident zeros and repeated copies.
The test functions are entire and rapidly decreasing in every fixed
horizontal strip, so this definition applies to complex ordinates. The
unitary functional equation pairs $\gamma_H$ with
$\overline{\gamma_H}$ for $H=F,G$. It follows that
$\mathcal R_{F,G}(T)$ is real. The evenness of $\psi$ also gives
$\mathcal R_{F,G}(T)=\mathcal R_{G,F}(T)$.

For an integrable function against $h^2$, use the notation
\begin{equation}\label{eq:eta}
 \av{v}=\frac1\kappa\int_\R h(t)^2v(t)\,dt,
 \qquad \eta_H(t)=\log Q_H+D_H\log\frac{|t|}{2\pi c}
 \quad(H=F,G).
\end{equation}
The logarithmic singularity at $t=0$ is integrable, including in the
product $\eta_F\eta_G$.

For use in the proofs, denote the absolute mixed prime-power convergence by
\begin{equation*}\tag{$\mathrm H_{F,G}^{\rm abs}$}\label{eq:Habs}
 \sum_{p,k\geq2}\frac{(k-1)(\log p)^2
              |w_{F,G}(p^k)|}{p^k}<\infty.
\end{equation*}
This convergence follows from Jiang's rational-field fixed-data
estimates \cite{Jiang}, after absorbing the logarithmic weight and
applying Cauchy--Schwarz to each cuspidal cross pair. Thus it is an
established input for the fixed products over $\Q$ considered here;
Proposition~\ref{prop:sufficient} gives the deduction in detail.
Set
\begin{equation}\label{eq:Aseries}
 A_{F,G}(s)=\sum_{p,k\geq1}
       \frac{(\log p)^2 w_{F,G}(p^k)}{p^{ks}}
       \qquad(\Rea s>1).
\end{equation}

\begin{theorem}[Secondary term]\label{thm:qualitative}
For the fixed products in \eqref{eq:products}, the limit
\begin{equation}\label{eq:C}
 C_{F,G}=\lim_{s\to1^+}
       \left(A_{F,G}(s)-\frac{r_{F,G}}{(s-1)^2}\right)
\end{equation}
exists, and
\begin{equation}\label{eq:main}
\begin{split}
 \mathcal R_{F,G}(T)={}&\frac{D_FD_G}{c^2}\phi(0)+r_{F,G}J_\phi\\
 &+\frac{\{D_F\av{\eta_G}+D_G\av{\eta_F}\}\phi(0)}{c^2L}\\
 &+\frac{\{\av{\eta_F\eta_G}+2\Rea C_{F,G}\}\phi(0)}{c^2L^2}
       +o(L^{-2}).
\end{split}
\end{equation}
\end{theorem}

For a PNT-based quantitative statement, let $\mathcal T\subset\R$
be the finite set consisting of $0$ and the norm-twist parameters
between constituents of $F$ and $G$. Their nonnegative integer
multiplicities are denoted by $r_\tau$, with $r_0=r_{F,G}$; zero
coefficients can be omitted. Harcos--Thorner \cite[Theorem 2.1]{HT}
gives, for every $A>0$,
\begin{equation}\label{eq:PNTrapid}
 \Psi_{F,G}(x):=\sum_{n\leq x}\Lambda(n)b_F(n)\overline{b_G(n)}
 =\sum_{\tau\in\mathcal T}r_\tau\frac{x^{1+i\tau}}{1+i\tau}
       +O_{F,G,A}\!\left(\frac{x}{(\log x)^A}\right).
\end{equation}
Here $\Lambda(n)$ is the von Mangoldt function. This estimate holds
for every fixed pair of products under consideration, with no
self-contragredience condition. Its constants may be ineffective.
Section~\ref{sec:established} explains the specialization to modulus
$1$ and the finite ramified correction.

We retain a separate notation for the stronger input
\begin{equation*}\tag{$\mathrm{PNT}_{F,G}$}\label{eq:PNT}
 \Psi_{F,G}(x)
 =\sum_{\tau\in\mathcal T}r_\tau\frac{x^{1+i\tau}}{1+i\tau}
       +O\!\left(xe^{-a_0\sqrt{\log x}}\right),
 \qquad a_0>0.
\end{equation*}
Theorem~2.3 of Liu--Ye \cite{LY07}, with the corrected proof and
ramified local treatment in the submitted corrigendum \cite{Ycorr},
proves this estimate whenever each cuspidal cross pair contains at
least one self-contragredient member. Thus both PNT statements are
established in their stated scopes. The general result
\eqref{eq:PNTrapid} does not assert the subexponential error in
\eqref{eq:PNT}.
The parameters correspond to pairs with
$\tau'\simeq\sigma\otimes|\det|^{i\tau}$, with multiplicities, where
$\sigma$ occurs in $F$ and $\tau'$ occurs in $G$. This convention agrees
with $a_{\sigma\otimes|\det|^{it}}(p^k)=p^{-ikt}a_\sigma(p^k)$.
The stronger mixed tail, proved in Proposition~\ref{prop:sufficient}, is
\begin{equation*}\tag{$\mathrm H_{F,G}^{\delta}$}\label{eq:Hdelta}
 \sum_{\substack{p^k>Y\\k\geq2}}
 \frac{(k-1)(\log p)^2|w_{F,G}(p^k)|}{p^k}
       \ll Y^{-\delta}\qquad(Y\geq2).
\end{equation*}
For every fixed pair of products over $\Q$, some $\delta>0$ is
available from Jiang's stronger rational-field fixed-data estimate
\eqref{eq:Jianginput}. Proposition~\ref{prop:sufficient} derives the
mixed tail by weighted Cauchy--Schwarz and gives the strict range
\eqref{eq:deltarange}. Thus the power saving is supplied by that
stronger estimate, rather than by qualitative Hypothesis H alone.
Finally, require a plateau:
\begin{equation}\label{eq:plateau}
 \phi(u)=\phi(0)\quad (|u|\leq a),\qquad 0<a<b.
\end{equation}

\begin{theorem}[Plateau refinements]\label{thm:quantitative}
Let $F,G$ be any fixed products in \eqref{eq:products}, and assume
the plateau condition \eqref{eq:plateau}. Let $\mathcal M_{F,G}(T)$
denote the four main terms displayed in \eqref{eq:main}. Then, for
every $N>0$,
\begin{equation}\label{eq:rapidremainder}
 \mathcal R_{F,G}(T)=\mathcal M_{F,G}(T)
 +\Theta_{F,G}(\lambda)+O_N(L^{-N})
 =\mathcal M_{F,G}(T)+O_N(L^{-N}),
\end{equation}
where
\begin{equation}\label{eq:theta}
 \Theta_{F,G}(\lambda)=2\Rea\sum_{\substack{\tau\in\mathcal T\\\tau\ne0}}
 r_\tau\left\{\int_0^\infty u\phi(u)e^{i\tau\lambda u}\,du
                 +\frac{\phi(0)}{\tau^2\lambda^2}\right\}.
\end{equation}
For every $N>0$, $\Theta_{F,G}(\lambda)=O_N(\lambda^{-N})$.
If the stronger Liu--Ye input \eqref{eq:PNT} holds, then, for some
$a_1>0$,
\begin{equation}\label{eq:quantremainder}
 \mathcal R_{F,G}(T)=\mathcal M_{F,G}(T)
 +\Theta_{F,G}(\lambda)+O\!\left(e^{-a_1\sqrt L}\right).
\end{equation}
Under that input, the twist term can be absorbed in a subexponential
remainder if $\mathcal T\setminus\{0\}$ is empty or if $\phi$ has
Gevrey order two. For a Gevrey-order-two plateau, the conclusion
$\mathcal R_{F,G}(T)=\mathcal M_{F,G}(T)+O(e^{-a_2\sqrt L})$
in fact holds for every fixed pair of products, by
Theorem~\ref{thm:localGevrey}, without \eqref{eq:PNT}.
Here the Gevrey condition means that there are constants $C_0,A_0>0$ such that
\begin{equation}\label{eq:gevrey}
 \|\phi^{(j)}\|_\infty\leq C_0A_0^j(j!)^2\qquad(j\geq0).
\end{equation}
The tail \eqref{eq:Hdelta} used in these assertions is automatic by
Proposition~\ref{prop:sufficient}.
\end{theorem}

\begin{remark}\label{rem:interpretation}
There is no arithmetic term of order $L^{-1}$ in \eqref{eq:main}.
That coefficient comes entirely from the smooth density. The first
arithmetic secondary coefficient is $2\Rea C_{F,G}$ at order $L^{-2}$.
With the established power tail, the general expansion in
Theorem~\ref{n:all-orders} improves the remainder for an arbitrary
even smooth pair test to $O(L^{-4})+O(T^{-\varepsilon})$.
For an arbitrary smooth plateau, \eqref{eq:rapidremainder} is
unconditional for the fixed products considered here. Rapid
inverse-power decay does not by itself imply a subexponential rate.
The stronger PNT with no nonzero twists, or the local Gevrey argument,
supplies that further refinement.
\end{remark}

\begin{corollary}[Cuspidal probe]\label{cor:probe}
Let $F$ be as in \eqref{eq:products}, and let $\pi$ be any fixed
unitary cuspidal representation over $\Q$. Put $G=L(s,\pi)$.
For every smooth plateau satisfying \eqref{eq:support},
\[
 \mathcal R_{F,G}(T)=\mathcal M_{F,G}(T)+O_N(L^{-N})
 \qquad\text{for every }N>0.
\]
For a Gevrey-order-two plateau, the error is
$O(e^{-a\sqrt{\log T}})$ for some $a>0$.
No self-contragredience, additional Hypothesis H, or temperedness
assumption is required. If $\pi$ is self-contragredient, Liu--Ye's
PNT additionally supplies \eqref{eq:quantremainder} for every smooth
plateau, with its nonzero norm-twist terms displayed explicitly.
\end{corollary}

\begin{corollary}[Self correlation]\label{cor:self}
For every fixed product $F$ in \eqref{eq:products},
\begin{equation}\label{eq:self}
\begin{split}
 \mathcal R_{F,F}(T)={}&\frac{D_F^2}{c^2}\phi(0)+B_FJ_\phi
 +\frac{2D_F\av{\eta_F}\phi(0)}{c^2L}\\
 &+\frac{\{\av{\eta_F^2}+2C_{F,F}\}\phi(0)}{c^2L^2}
 +o(L^{-2}).
\end{split}
\end{equation}
The constant $C_{F,F}$ is real. For every smooth plateau, the
remainder can be replaced by $O_N(L^{-N})$ for every $N>0$.
For a Gevrey-order-two plateau, it can be replaced by
$O(e^{-a\sqrt{\log T}})$ for some $a>0$. If every constituent of
$F$ is self-contragredient, the latter remainder holds for every
smooth plateau.
\end{corollary}

\section{The arithmetic finite part}

We first record bounds that require neither Hypothesis H nor a
quantitative prime number theorem.

\begin{lemma}\label{lem:bounds}
For a fixed product $H$ as in \eqref{eq:products},
\begin{equation}\label{eq:diagonalbound}
 \sum_{n\leq x}\Lambda(n)|b_H(n)|^2\ll x,\qquad
 \sum_{n\leq x}\Lambda(n)|b_H(n)|\ll x.
\end{equation}
For two such products,
\begin{equation}\label{eq:coeffbound}
 \sum_{n\leq x}\frac{\Lambda(n)^2|b_F(n)b_G(n)|}{n}
       \ll (\log(2x))^2.
\end{equation}
Moreover, for $\lambda\geq1$,
\begin{equation}\label{eq:globalBbound}
 \sum_{n\geq2}\frac{(\log n)\Lambda(n)|b_F(n)b_G(n)|}
                         {n^{1+1/\lambda}}\ll\lambda^2.
\end{equation}
\end{lemma}

\begin{proof}
For a single unitary cuspidal $\sigma$, the series
$\sum\Lambda(n)|a_\sigma(n)|^2n^{-s}$ differs from
$-L'/L(s,\sigma\times\wt\sigma)$ by a finite ramified correction
holomorphic in a half-plane $\Rea s>\theta$ with $\theta<1$.
Indeed, at each ramified prime the standard roots have absolute value
less than $p^{1/2}$ and the genuine Rankin--Selberg roots have absolute
value less than $p$; see \cite{BTZ,Ycorr}. The series has nonnegative
coefficients. The simple pole at $1$ and nonvanishing on $\Rea s=1$
therefore allow the Wiener--Ikehara theorem to give
\[
 \sum_{n\leq x}\Lambda(n)|a_\sigma(n)|^2\sim x.
\]
The required absolute convergence for $\Rea s>1$ follows first from
the diagonal Rankin--Selberg Euler product and the finite local
correction; thus this argument does not assume the bound it proves.

For the product $H$, apply
\[
 \left|\sum_\sigma t_\sigma a_\sigma(n)\right|^2
 \leq\left(\sum_\sigma t_\sigma\right)
       \sum_\sigma t_\sigma|a_\sigma(n)|^2.
\]
This proves the first assertion of \eqref{eq:diagonalbound}.
The second follows by Cauchy--Schwarz and
$\sum_{n\leq x}\Lambda(n)\ll x$. Another application of
Cauchy--Schwarz gives
$\sum_{n\leq x}\Lambda(n)|b_F(n)b_G(n)|\ll x$.
Partial summation, together with $\Lambda(n)\leq\log n$, proves
\eqref{eq:coeffbound} and \eqref{eq:globalBbound}.

The constituentwise argument is useful: applying a one-pole Tauberian
argument directly to the aggregate product would overlook possible
poles at $1+i\tau$ from distinct norm-related constituents.
\end{proof}

Define, for $\Rea s>1$,
\begin{align}
 K_{F,G}(s)&=\sum_{p,k\geq1}\frac{(\log p)w_{F,G}(p^k)}{p^{ks}},
       \label{eq:K}\\
 H_{F,G}(s)&=\sum_{p,k\geq2}
       \frac{(k-1)(\log p)^2w_{F,G}(p^k)}{p^{ks}}.
       \label{eq:Hseries}
\end{align}
Then the exact identity is
\begin{equation}\label{eq:Aderivative}
 A_{F,G}(s)=-K_{F,G}'(s)-H_{F,G}(s).
\end{equation}
The factor $k-1$ corrects the factor $k$ obtained by differentiating
$p^{-ks}$. It cannot be omitted in a secondary-term calculation.

Let $S$ contain the primes ramified in any constituent of $F$ or $G$,
and put
\[
 \mathcal L^S_{F,G}(s)=
 \prod_{\sigma,\tau\in\mathcal S}
 L^S(s,\sigma\times\wt\tau)^{r_\sigma s_\tau}.
\]
The superscript means that the local factors at $S$ have been removed.
At unramified primes the tensor-product identity for the roots gives
\begin{equation}\label{eq:badprimes}
 K_{F,G}(s)=-\frac{d}{ds}\log\mathcal L^S_{F,G}(s)
 +\sum_{p\in S,k\geq1}\frac{(\log p)w_{F,G}(p^k)}{p^{ks}}.
\end{equation}
The last series is holomorphic in a neighborhood of $\Rea s\geq1$.
The order of the pole of $\mathcal L^S_{F,G}$ at $1$ is $r_{F,G}$.
Consequently, with $r=r_{F,G}$,
\begin{equation}\label{eq:Laurent}
 K_{F,G}(1+z)=\frac r z+\beta_0+\beta_1z+O(z^2)
\end{equation}
in a fixed punctured disk. Under \eqref{eq:Habs},
$H_{F,G}(1+z)\to H_{F,G}(1)$ as $z\to0^+$, and
\begin{equation}\label{eq:Cbeta}
 C_{F,G}=-\beta_1-H_{F,G}(1).
\end{equation}
This proves existence in \eqref{eq:C}. It also shows that the expansion
of $A_{F,G}(1+z)$ as real $z\to0^+$ contains no $z^{-1}$ term.

For clarity, an alternative expression displaying all ramified terms is
\begin{equation}\label{eq:Cexplicit}
\begin{split}
 C_{F,G}={}&
 \left.\frac{d^2}{ds^2}
       \log\big((s-1)^r\mathcal L^S_{F,G}(s)\big)\right|_{s=1}\\
 &-\sum_{\substack{p\notin S\\k\geq2}}
       \frac{(k-1)(\log p)^2w_{F,G}(p^k)}{p^k}
 +\sum_{\substack{p\in S\\k\geq1}}
       \frac{(\log p)^2w_{F,G}(p^k)}{p^k}.
\end{split}
\end{equation}
The logarithm denotes any local holomorphic branch of the nonvanishing
regularized function. Its second derivative is independent of the
branch. Formula \eqref{eq:Cexplicit} does not identify ramified
Rankin--Selberg roots with products of standard local roots.

\begin{lemma}[Qualitative smoothed sum]\label{lem:smoothed}
Assume \eqref{eq:Habs}. For $W\in C_c^\infty(\R)$,
\begin{equation}\label{eq:smoothed}
 \sum_{p,k\geq1}\frac{(\log p)^2w_{F,G}(p^k)}{p^k}
 W\!\left(\frac{k\log p}{\lambda}\right)
 =r\lambda^2\int_0^\infty uW(u)\,du+C_{F,G}W(0)+o(1).
\end{equation}
\end{lemma}

\begin{proof}
Write $K=K_{F,G}$ and put $f(u)=e^uW(u)$,
$\widehat f(t)=\int_\R f(u)e^{-itu}\,du$. Fourier inversion gives
\begin{equation}\label{eq:Fourierlocal}
 \sum_{p,k\geq1}\frac{k(\log p)^2w_{F,G}(p^k)}{p^k}
 W\!\left(\frac{k\log p}{\lambda}\right)
 =\frac1{2\pi}\int_\R\widehat f(t)
       [-K'(1+(1-it)/\lambda)]\,dt.
\end{equation}
The interchange is absolutely justified by
\eqref{eq:globalBbound}. That bound is uniform in the imaginary part
and bounds the bracket in absolute value by $O(\lambda^2)$.

Choose a fixed disk about $1$ in which \eqref{eq:Laurent} holds and
choose $\varepsilon_0>0$ sufficiently small. On
$|t|\leq\varepsilon_0\lambda$,
\[
 -K'(1+(1-it)/\lambda)
 =\frac{r\lambda^2}{(1-it)^2}-\beta_1
       +O((1+|t|)/\lambda).
\]
The complementary part of the integral is smaller than every prescribed
inverse power of $\lambda$, after using the Schwartz decay of
$\widehat f$ and the global $O(\lambda^2)$ bound. Fourier inversion
and the Laplace representation of $(1-it)^{-2}$ give
\[
 \frac1{2\pi}\int_\R\frac{\widehat f(t)}{(1-it)^2}\,dt
 =\int_0^\infty u e^{-u}f(u)\,du
 =\int_0^\infty uW(u)\,du.
\]
Also $(2\pi)^{-1}\int\widehat f(t)\,dt=W(0)$. Thus the left side of
\eqref{eq:Fourierlocal} equals
\[
 r\lambda^2\int_0^\infty uW(u)\,du-\beta_1W(0)+O(\lambda^{-1}).
\]
Subtract the sum with coefficient $k-1$. By absolute convergence in
\eqref{eq:Habs} and dominated convergence, that sum tends to
$H_{F,G}(1)W(0)$. Equation \eqref{eq:Cbeta} proves
\eqref{eq:smoothed}. Only a local Laurent expansion was used; no
global zero-free region or contour shift is needed for this lemma.
\end{proof}

\section{Reduction of the zero sum}

\begin{proposition}\label{prop:reduction}
Under \eqref{eq:products}--\eqref{eq:support}, without
\eqref{eq:Habs} or \eqref{eq:PNT}, there is an $\varepsilon>0$ such that
\begin{equation}\label{eq:reduction}
 \mathcal R_{F,G}(T)=
 \frac{\phi(0)}{\lambda^2}
       \av{(D_FL+\eta_F)(D_GL+\eta_G)}
 +\frac2{\lambda^2}\Rea S_{F,G}(\lambda)+O(T^{-\varepsilon}),
\end{equation}
where
\begin{equation}\label{eq:S}
 S_{F,G}(\lambda)=\sum_{p,k\geq1}
       \frac{(\log p)^2w_{F,G}(p^k)}{p^k}
       \phi\!\left(\frac{k\log p}{\lambda}\right).
\end{equation}
One may take any fixed $0<\varepsilon<(1-cb)/2$.
\end{proposition}

\begin{proof}
We give the terms in the explicit formula separately. Define the exact
archimedean density by
\begin{equation}\label{eq:Omega}
\begin{split}
 \Omega_F(t)=\log Q_F+\sum_{j=1}^{D_F}\bigg\{
 &\frac{\Gamma_\R'}{\Gamma_\R}
             (\tfrac12+\mu_{F,j}+it)\\
 +&\frac{\Gamma_\R'}{\Gamma_\R}
             (\tfrac12+\overline{\mu_{F,j}}-it)\bigg\}.
\end{split}
\end{equation}
It is real and smooth on $\R$. The standard archimedean parameters of
a unitary cuspidal representation satisfy $\Rea\mu_{F,j}>-1/2$,
so there are no gamma poles on this line. Put
\[
 q_{F,U}(x)=\frac1{2\pi}\int_\R
                   h(t)\Omega_F(Ut)e^{-ixt}\,dt,
 \qquad Z_F(u)=\sum_{\gamma_F}h(\gamma_F/U)e^{-i\lambda u\gamma_F}.
\]
The explicit formula, with the Fourier conventions
\eqref{eq:h}, gives
\begin{equation}\label{eq:EF}
\begin{split}
 Z_F(u)={}&Uq_{F,U}(U\lambda u)+P_F(u)\\
 &-U\sum_{n\geq2}\frac{\Lambda(n)b_F(n)}{\sqrt n}
                          g(U\lambda u+U\log n)\\
 &-U\sum_{n\geq2}\frac{\Lambda(n)\overline{b_F(n)}}{\sqrt n}
                          g(U\lambda u-U\log n).
\end{split}
\end{equation}
Here $P_F$ is the finite sum from the poles of the completed function,
with their multiplicities. This is the scaled form of the explicit
formula used in \cite{RS}. Any such poles come from $\GL_1$ norm
characters; their real parts are $0$ and $1$, with fixed possible
vertical shifts. Hence, uniformly for $|u|\leq b$,
\begin{equation}\label{eq:Pbound}
 |P_F(u)|\ll e^{\lambda|u|/2}.
\end{equation}
The same notation and estimates apply to $G$.

By evenness of $\phi$, the unnormalized sum in
\eqref{eq:correlation} is
\begin{equation}\label{eq:rawintegral}
 \int_\R\phi(u)Z_F(u)Z_G(-u)\,du.
\end{equation}
Every arithmetic packet occurring in this integral has
\begin{equation}\label{eq:X}
 n\leq X:=\exp(\lambda b+A/U)=T^{cb+o(1)}.
\end{equation}
Two packets centered at the same side of the origin can overlap only if
$U|\log(n/m)|\leq2A$. For distinct integers $n,m\leq X$,
\[
 U|\log(n/m)|\gg U/X\longrightarrow\infty.
\]
Thus unequal arithmetic packets are exactly disjoint for sufficiently
large $T$. Packets centered on opposite sides cannot overlap either,
since $U\log(nm)\geq U\log4>2A$.

For the matching orientation and $n=m$, the packet integral, including
the factor $U^2$, is
\[
 \frac U\lambda\int_\R
 \phi\!\left(\frac{\log n}{\lambda}+\frac v{U\lambda}\right)
                          g(v)g(-v)\,dv.
\]
Replacing $\phi$ by $\phi(\log n/\lambda)$ costs at most
$O(\lambda^{-2})$ times the coefficient sum in
\eqref{eq:coeffbound}, and therefore $O(1)$ in
\eqref{eq:rawintegral}. Since $\int g(v)g(-v)\,dv=\kappa/(2\pi)$,
the two matching orientations yield, after normalization,
\begin{equation}\label{eq:arithmeticpart}
 \frac2{\lambda^2}\Rea S_{F,G}(\lambda)+O((TL)^{-1}).
\end{equation}

We next evaluate the product of the two archimedean terms. The gamma
estimates and rapid decay of $h$ imply
\[
 \|h(t)\Omega_F(Ut)\|_1+
 \left\|\frac d{dt}\{h(t)\Omega_F(Ut)\}\right\|_1\ll L.
\]
For example, $\Omega_F(v)=D_F\log(|v|/(2\pi))+\log Q_F+O(|v|^{-1})$
for large $|v|$, and $\Omega_F'(v)=O((1+|v|)^{-1})$.
Integration by parts gives
\begin{equation}\label{eq:qbound}
 |q_{F,U}(x)|\ll\frac{L}{1+|x|}.
\end{equation}
The raw archimedean contribution is
\[
 \frac U\lambda\int_\R
       \phi\!\left(\frac{x}{U\lambda}\right)
                    q_{F,U}(x)q_{G,U}(-x)\,dx.
\]
Since $\phi$ is even,
\[
 \left|\phi\!\left(\frac{x}{U\lambda}\right)-\phi(0)\right|
 \ll\min\left\{\frac{x^2}{(U\lambda)^2},1\right\}.
\]
Together with \eqref{eq:qbound}, this bounds the cost of replacing
$\phi$ by $\phi(0)$ by $O(UL^2/(\lambda U\lambda))=O(1)$.
Parseval's identity now gives the normalized contribution
\begin{equation}\label{eq:exactdensity}
 \frac{\phi(0)}{\lambda^2}\av{\Omega_F(Ut)\Omega_G(Ut)}
       +O((TL)^{-1}).
\end{equation}
For $|Ut|\geq1$, Stirling's formula gives
\[
 \Omega_H(Ut)=D_HL+\eta_H(t)+O((U|t|)^{-1}),\qquad H=F,G.
\]
On $|t|<1/U$, use smoothness of $\Omega_H$ and integrate the
logarithmic terms directly. On the complementary interval integrate
the displayed remainder using the decay of $h$. This proves
\begin{equation}\label{eq:densityapprox}
 \av{\Omega_F(Ut)\Omega_G(Ut)}
 =\av{(D_FL+\eta_F)(D_GL+\eta_G)}+O(L^2/U).
\end{equation}

It remains to bound all cross terms. On an arithmetic packet,
$|U\lambda u|\asymp U\log n$, and its width is $O((U\lambda)^{-1})$.
Using \eqref{eq:qbound}, \eqref{eq:diagonalbound}, and partial
summation, a gamma--arithmetic term is $O(L\sqrt X)$ in the raw
integral (a sharper bound is unnecessary). A pole--pole term is
$O(X)$ by \eqref{eq:Pbound}. On the packet for $n$, that same bound
is $O(\sqrt n)$; therefore a pole--arithmetic term is bounded by
\[
 \frac1\lambda\sum_{n\leq X}\Lambda(n)|b_H(n)|\ll X/\lambda.
\]
This averaged bound avoids a loss from a pointwise estimate for the
coefficients. Finally, by \eqref{eq:qbound}, a pole--gamma term is
\[
 \ll\frac L\lambda\int_{|x|\leq U\lambda b}
             \frac{e^{|x|/(2U)}}{1+|x|}\,dx
 \ll L\sqrt X.
\]
After division by $TL$, all these terms are
$O(T^{-\varepsilon})$ for any
$0<\varepsilon<(1-cb)/2$. Combining
\eqref{eq:arithmeticpart}--\eqref{eq:densityapprox} proves the
proposition.
\end{proof}

\begin{proof}[Proof of Theorem~\ref{thm:qualitative}]
Apply Lemma~\ref{lem:smoothed} with $W=\phi$ in
Proposition~\ref{prop:reduction}. Since $\phi$ is real and even,
$2\int_0^\infty u\phi(u)\,du=J_\phi$.
Expand the density product and use $\lambda=cL$. This gives
\eqref{eq:main}.
\end{proof}

\section{Plateau refinements and norm twists}

The plateau condition allows the arithmetic remainder to be estimated
using only the tail of a Mertens remainder. We give the calculation with
norm twists included, because its constant differs from the untwisted
one.

\begin{lemma}\label{lem:quantarith}
Let $F,G$ be any fixed products in \eqref{eq:products}, and let
$W\in C_c^\infty(\R)$ be constant with value $W(0)$ on $[0,a]$ for
some $a>0$. Then, for every $N>0$,
\begin{equation}\label{eq:quantsum}
\begin{split}
 &\sum_{p,k\geq1}\frac{(\log p)^2w_{F,G}(p^k)}{p^k}
                   W\!\left(\frac{k\log p}{\lambda}\right)\\
 &\quad=r_0\lambda^2\int_0^\infty uW(u)\,du+C_{F,G}W(0)\\
 &\qquad\quad+\sum_{\substack{\tau\in\mathcal T\\\tau\ne0}}
 r_\tau\left\{\lambda^2\int_0^\infty uW(u)e^{i\tau\lambda u}\,du
                       +\frac{W(0)}{\tau^2}\right\}
 +O_N(\lambda^{-N}).
\end{split}
\end{equation}
If \eqref{eq:PNT} holds, the final error can be replaced by
$O(e^{-a_3\sqrt\lambda})$ for some $a_3>0$.
\end{lemma}

\begin{proof}
For $y\geq0$, define the right-continuous function
\[
 M(y)=\sum_{n\leq e^y}
                 \frac{\Lambda(n)b_F(n)\overline{b_G(n)}}n.
\]
Partial summation in \eqref{eq:PNTrapid}, with its logarithmic
exponent chosen sufficiently large, gives
\begin{equation}\label{eq:Mertens}
 M(y)=r_0y+\sum_{\substack{\tau\in\mathcal T\\\tau\ne0}}
             \frac{r_\tau e^{i\tau y}}{i\tau}+\beta_0+\epsilon(y),
 \qquad \epsilon(y)\ll_A(1+y)^{-A}
\end{equation}
for every $A>0$. Indeed, the error in the weighted partial sum is
bounded by the endpoint error divided by $e^y$ and the tail integral
of the PNT error divided by $x^2$; a saving $(\log x)^{-A-1}$
therefore suffices. Under \eqref{eq:PNT}, the same calculation gives
$\epsilon(y)\ll e^{-a_4\sqrt y}$ for some $a_4>0$, after decreasing
the exponent constant to absorb a polynomial factor. Both estimates
have a harmless bounded interpretation at small $y$.
To identify the constant with \eqref{eq:Laurent}, observe for real
$z>0$ that
\[
 K_{F,G}(1+z)=z\int_0^\infty M(y)e^{-zy}\,dy.
\]
The nonzero-twist term contributes
\[
 \frac{r_\tau z}{i\tau(z-i\tau)}
      =\frac{r_\tau z}{\tau^2}+O(z^2).
\]
The constant term in this Laplace expansion is precisely $\beta_0$,
and its linear coefficient is
\begin{equation}\label{eq:beta1twist}
 \beta_1=\sum_{\substack{\tau\in\mathcal T\\\tau\ne0}}
                      \frac{r_\tau}{\tau^2}
                   +\int_0^\infty\epsilon(y)\,dy.
\end{equation}
The integral converges absolutely, since we may take $A>1$ in
\eqref{eq:Mertens}. In the absence of nonzero twists,
\eqref{eq:beta1twist} reduces to
$\beta_1=\int_0^\infty\epsilon(y)\,dy$.

The sum with the additional weight $k$ is the Stieltjes integral
\[
 \int_0^\infty yW(y/\lambda)\,dM(y).
\]
The term $r_0y$ gives
$r_0\lambda^2\int_0^\infty uW(u)\,du$, and each exponential in
\eqref{eq:Mertens} gives
$r_\tau\lambda^2\int_0^\infty uW(u)e^{i\tau\lambda u}\,du$.
Integration by parts for the remainder gives
\[
 -\int_0^\infty\epsilon(y)
       \{W(y/\lambda)+(y/\lambda)W'(y/\lambda)\}\,dy.
\]
The constant $\beta_0$ makes no contribution since $yW(y/\lambda)$
vanishes at both endpoints. The expression in braces equals $W(0)$
on $0\leq y\leq a\lambda$ and is bounded uniformly elsewhere.
It follows that the last integral equals
\[
 -W(0)\int_0^\infty\epsilon(y)\,dy
       +O_N(\lambda^{-N})
\]
for every $N>0$. More explicitly, the error is bounded by a constant
times $\int_{a\lambda}^\infty|\epsilon(y)|\,dy$, which is
$O_A(\lambda^{1-A})$; take $A>N+1$.
Under \eqref{eq:PNT}, the same tail integral is
$O(e^{-a_5\sqrt\lambda})$, after decreasing the positive constant
in the exponent. Use \eqref{eq:beta1twist} to rewrite its constant.

Finally, subtract the sum with weight $k-1$. Its difference from
$H_{F,G}(1)W(0)$ is supported on $p^k>e^{a\lambda}$ and is
$O(e^{-\delta a\lambda})$ by \eqref{eq:Hdelta}. The identity
\eqref{eq:Cbeta} now gives \eqref{eq:quantsum}.
\end{proof}

\begin{lemma}\label{lem:Gevrey}
If $\phi$ satisfies \eqref{eq:plateau}, then for each fixed
$\tau\ne0$,
\begin{equation}\label{eq:twistIBP}
 \int_0^\infty u\phi(u)e^{i\tau\lambda u}\,du
       +\frac{\phi(0)}{\tau^2\lambda^2}
 =-\frac1{\tau^2\lambda^2}
    \int_0^\infty\{2\phi'(u)+u\phi''(u)\}e^{i\tau\lambda u}\,du.
\end{equation}
It is $O_N(\lambda^{-N})$ for every $N$. Under
\eqref{eq:gevrey}, it is $O(e^{-a_6\sqrt\lambda})$ for some $a_6>0$.
For every $0<a<b$ there exists an even Gevrey-order-two test that is
one on $[-a,a]$ and supported in $[-b,b]$.
\end{lemma}

\begin{proof}
Two integrations by parts give \eqref{eq:twistIBP}. The function
$2\phi'+u\phi''$ on the positive half-line is smooth, compactly
supported, and vanishes near zero. It can therefore be integrated by
parts arbitrarily many further times without boundary terms.

Under \eqref{eq:gevrey}, the $N$th derivative of this function has
$L^1$ norm at most $CA^N(N!)^2$ after changing $C,A$. Hence its
Fourier integral at frequency $\omega$ is bounded by
$CA^N(N!)^2|\omega|^{-N}$. Taking
$N\asymp\sqrt{|\omega|}$ proves exponential decay in
$\sqrt{|\omega|}$, and $\omega=\tau\lambda$ proves the assertion.

For existence, let $\beta(x)=e^{-1/x}$ for $x>0$ and $\beta(x)=0$
for $x\leq0$. On bounded intervals this function is of Gevrey order
two: Cauchy's formula on a circle $|z-x|=\alpha x$, with a fixed
small $\alpha>0$, bounds its $j$th derivative by
\[
 j!(\alpha x)^{-j}e^{-c_0/x}\leq CA^j(j!)^2.
\]
Normalize $\rho(v)=\beta(1-v)\beta(1+v)$ to have integral one.
This is an even Gevrey-order-two bump supported in $[-1,1]$.
Put $\epsilon=(b-a)/2$, $d=(a+b)/2$, and
$\rho_\epsilon(v)=\epsilon^{-1}\rho(v/\epsilon)$. Then
\[
 \phi=\mathbf1_{[-d,d]}*\rho_\epsilon
\]
has the claimed properties and satisfies \eqref{eq:gevrey}.
\end{proof}

\begin{proof}[Proof of Theorem~\ref{thm:quantitative}]
Insert \eqref{eq:quantsum} with $W=\phi$ into
Proposition~\ref{prop:reduction}. The nonzero-twist terms are exactly
\eqref{eq:theta}. The general remainder in \eqref{eq:quantsum},
multiplied by $2\lambda^{-2}$, and the power saving in $T$ give
\eqref{eq:rapidremainder}, since $\lambda=c\log T$ and the chosen
inverse-power order is arbitrary. Under \eqref{eq:PNT}, the
subexponential arithmetic error and the power saving in $T$ are
absorbed in $O(e^{-a_1\sqrt{\log T}})$, proving
\eqref{eq:quantremainder}. Lemma~\ref{lem:Gevrey} treats the twist
term. The assertion for arbitrary products and a Gevrey plateau is
proved independently in Theorem~\ref{thm:localGevrey}.
\end{proof}

\section{Established arithmetic inputs and a stronger pair refinement}
\label{sec:established}

The symbols $\mathrm H^{\rm abs}$, $\mathrm H^\delta$, and
$\mathrm{PNT}$ distinguish the analytic inputs used in the preceding
proofs. The first two and the general PNT \eqref{eq:PNTrapid} are
established for all the fixed products considered here; the stronger
subexponential PNT \eqref{eq:PNT} is established under the Liu--Ye
hypotheses. We now give their precise sources. Throughout this section the
representations are fixed, and all the constituent representations are
over $\Q$.

\subsection{The quantitative consequence of Jiang's result}

For a cuspidal representation $\sigma$ of degree $d_\sigma$, put
$\epsilon_\sigma=(d_\sigma^2+1)^{-1}$. We use the rational-field case of
Jiang's Theorem~1.3 \cite{Jiang}. Its two assertions, together with
inequality (3.3) of that paper, give, for every
$0<\eta<\epsilon_\sigma$,
\begin{equation}\label{eq:Jianginput}
 \sum_p\sum_{k\ge2}
 \frac{|a_\sigma(p^k)|^2+
       \lambda_{\sigma\times\widetilde\sigma}(p^k)}
      {p^{k(1-\eta)}}<\infty.
\end{equation}
Here $\lambda_{\sigma\times\widetilde\sigma}$ denotes the ordinary
Dirichlet coefficients of the Rankin--Selberg $L$-function, not the
coefficients of its logarithmic derivative. Positivity turns the
convergent Euler products in that theorem into convergent sums. The
second assertion uses
$|a_\sigma(p^k)|^2\le a_{\sigma\times\widetilde\sigma}(p^k)$.

Only fixed-representation convergence in \eqref{eq:Jianginput} is used
here. In particular, we do not require the conductor-uniform iteration
or a number-field extension of the power sieve. For this purpose one
can use the rational-field estimate (5.6) in \cite{Jiang} for the
finitely many small exponents, followed by the local root bound for
large exponents. In (5.3) of that preprint, the displayed bound for the
$L$-value approaching its pole must retain its dependence on
$s-1$: for fixed $\sigma$ the Rankin argument gives
$O_\sigma(X\log X)$. This logarithm is absorbed by the arbitrarily
small powers already present in (5.6). Thus the fixed-data input used
here is unaffected. We cite the version reviewed here, and do not rely
on its general-number-field character construction.

\begin{proposition}[The mixed power tail]\label{prop:sufficient}
Let $M_F$ and $M_G$ be the largest constituent degrees in $F$ and $G$.
Then \eqref{eq:Hdelta}, and hence \eqref{eq:Habs}, holds for every
\begin{equation}\label{eq:deltarange}
 0<\delta<\frac12\left\{\frac1{M_F^2+1}+\frac1{M_G^2+1}\right\}.
\end{equation}
The implied constant may depend on the fixed products and on $\delta$.
\end{proposition}

\begin{proof}
Fix $0<\eta_\sigma<\epsilon_\sigma$ and apply
\eqref{eq:Jianginput} with a slightly larger exponent gain.
Since
\[
 (k-1)(\log p)^2\ll (\log(p^k))^3\ll_\varepsilon(p^k)^\varepsilon,
\]
the spare gain absorbs the weight and proves
\begin{equation}\label{eq:weightedJiang}
 M_\sigma(\eta_\sigma):=
 \sum_{p,k\ge2}\frac{(k-1)(\log p)^2|a_\sigma(p^k)|^2}
                         {p^{k(1-\eta_\sigma)}}<\infty.
\end{equation}
For two constituents, put $\delta=(\eta_\sigma+\eta_\tau)/2$.
Weighted Cauchy--Schwarz gives
\begin{align*}
 &\sum_{p,k\ge2}
 \frac{(k-1)(\log p)^2|a_\sigma(p^k)a_\tau(p^k)|}
                         {p^{k(1-\delta)}}\\
 &\qquad\le M_\sigma(\eta_\sigma)^{1/2}
               M_\tau(\eta_\tau)^{1/2}<\infty.
\end{align*}
On $p^k>Y$, insert $(p^k)^{-\delta}\le Y^{-\delta}$.
This proves the required tail for this pair. Expand the finite sums
defining $b_F,b_G$, use the triangle inequality, and choose the gains
uniformly over the finitely many cross pairs. This gives exactly the
strict range \eqref{eq:deltarange}. No assertion at its endpoint is made.
\end{proof}

Consequently, the qualitative secondary-term theorem is a theorem for
the stated fixed automorphic products over $\Q$, using the cited
rational-field result of Jiang. Hypothesis H is not an additional
conjectural assumption. The quantitative tail, which is stronger than
ordinary Hypothesis H, is also available. The ordinary-coefficient part
of \eqref{eq:Jianginput} will have a separate role in the higher-level
calculation: it controls fourth and higher powers of prime coefficients.

\subsection{The prime number theorems of Harcos--Thorner and Liu--Ye}

Apply Harcos--Thorner \cite[Theorem 2.1]{HT} over $\Q$, with modulus
$1$, to each pair $(\sigma,\widetilde\tau)$ of cuspidal
constituents of $F$ and $G$. This gives the genuine Rankin--Selberg
logarithmic-derivative coefficient sum with error
$O_{\sigma,\tau,A}(x(\log x)^{-A})$ for every $A>0$. No member of
the pair is required to be self-contragredient. At an unramified
prime these coefficients are
$(\log p)a_\sigma(p^k)\overline{a_\tau(p^k)}$.
At each of the finitely many remaining primes, the genuine
Rankin--Selberg roots and the products of standard roots have
modulus strictly less than $p$; see \cite[(3.9)--(3.10)]{HT} and the
local discussion in Section~3. Their partial-sum discrepancy is
therefore $O_{\sigma,\tau}(x^\theta\log x)$ for some fixed
$\theta<1$. This is absorbed in every stated logarithmic saving.
Summing the finitely many pairs proves \eqref{eq:PNTrapid}.
The multiplicity of a twist, in the convention of Section~2, is
\begin{equation}\label{eq:PNTrtau}
 r_t=\sum_{\substack{\sigma,\tau\in\mathcal S\\
               \tau\simeq\sigma\otimes|\det|^{it}}}
                  r_\sigma s_\tau.
\end{equation}
The constants in this general PNT may be ineffective. This causes
no difficulty because all representations and twists are fixed.

Theorem~2.3 of Liu--Ye \cite{LY07}, printed p.~485, with the
corrected proof and ramified local terms in \cite{Ycorr}, gives the
stronger subexponential estimate \eqref{eq:PNT} whenever at least
one member of each cuspidal cross pair is self-contragredient.
In particular, it suffices that every constituent of one product
is self-contragredient. The possible norm-twist terms are retained
in that theorem. The exceptional-zero issue is absorbed in the
constants for these fixed representations; no uniformity in
conductors or a moving twist is asserted. The standard zero-free
input in this setting is also discussed in Brumley's appendix to
\cite{Humphries}. Thus Harcos--Thorner removes the
self-contragredience restriction for the general PNT and its
rapid inverse-logarithmic consequences, while Liu--Ye supplies
the stronger error in its stated scope.

For clarity, these results also imply the weighted prime form of
Selberg orthogonality used in small-support correlation arguments.
Put
\[
 Q_{F,G}(x)=\sum_{\substack{p^k\le x\\k\ge2}}
       \frac{(\log p)w_{F,G}(p^k)}{p^k}.
\]
The established tail \eqref{eq:Hdelta}, together with
$\log p\le (k-1)(\log p)^2/\log2$, gives
$Q_{F,G}(x)=Q_{F,G}(\infty)+O(x^{-\delta})$.
Subtracting this from \eqref{eq:Mertens}, with $y=\log x$, yields
\begin{equation}\label{eq:primeMertens}
 \sum_{p\le x}\frac{(\log p)b_F(p)\overline{b_G(p)}}p
 =r_0\log x+\sum_{\substack{t\in\mathcal T\\t\ne0}}
      \frac{r_tx^{it}}{it}+C^{\rm pr}_{F,G}
      +O_N((\log x)^{-N})
\end{equation}
for every $N>0$, where $C^{\rm pr}_{F,G}=\beta_0-Q_{F,G}(\infty)$.
Under Liu--Ye's hypotheses the final error is
$O(e^{-a\sqrt{\log x}})$ for some $a>0$.
In particular, for a single cuspidal pair unrelated by any imaginary
norm twist, all $r_t$ vanish and the weighted prime sum is bounded.
This implication uses Jiang's quantitative prime-power tail in
addition to the respective PNT; qualitative Hypothesis H alone
would give convergence without the stated tail rate. This supplies
the prime orthogonality input in the corrected superposition results
of \cite{LY02,YForumCorr}.

\begin{proof}[Proof of Corollary~\ref{cor:probe}]
For any fixed cuspidal probe, \eqref{eq:PNTrapid} and
Proposition~\ref{prop:sufficient} apply to every cross pair. The
first assertion follows from Theorem~\ref{thm:quantitative}; the
Gevrey assertion follows from Theorem~\ref{thm:localGevrey}.
If the probe is self-contragredient, every cross pair also meets
the Liu--Ye hypothesis, giving \eqref{eq:quantremainder}.
\end{proof}

\begin{proof}[Proof of Corollary~\ref{cor:self}]
Put $G=F$ in Theorem~\ref{thm:qualitative}. All coefficients of
$A_{F,F}$ are real, so $C_{F,F}$ is real. The general smooth-plateau
assertion follows from \eqref{eq:rapidremainder}, and the Gevrey
assertion follows from Theorem~\ref{thm:localGevrey}. If every
constituent is self-contragredient, apply the Liu--Ye PNT to every
pair. Two such constituents cannot differ by a nonzero imaginary
norm twist: compare their central characters on the positive real
ideles. Thus the twist term is absent, and any smooth plateau gives
the asserted subexponential refinement.
\end{proof}

\subsection{Local Gevrey localization: the PNT is unnecessary here}

For the narrower class of Gevrey-order-two plateau tests, there is a
stronger statement. A local argument near $s=1$ replaces the global
prime number theorem entirely. This also handles nonzero norm twists.

\begin{theorem}[Pair refinement for arbitrary fixed products]
\label{thm:localGevrey}
Let $F,G$ be any fixed products in \eqref{eq:products}. If $\phi$ is an
even Gevrey-order-two plateau satisfying \eqref{eq:support}, then
\eqref{eq:main} holds with error
\begin{equation}\label{eq:unconditionalpairerror}
 O\!\left(e^{-a\sqrt{\log T}}\right)
\end{equation}
for some $a>0$. No self-contragredience, temperedness, or separate
assumption \eqref{eq:PNT} is required. The constants depend on the fixed
representations, including their norm-twist separations.
\end{theorem}

\begin{proof}
Proposition~\ref{prop:sufficient} makes $H_{F,G}(1+z)$ holomorphic in
a neighborhood of zero. The Rankin--Selberg factorization of $K$ and
nonvanishing at $1$, after removing its pole, therefore give a convergent
local expansion
\begin{equation}\label{eq:fullAjet}
 A_{F,G}(1+z)=\frac{r}{z^2}+\sum_{j\ge0}C_{F,G,j}z^j,
 \qquad C_{F,G,0}=C_{F,G},\qquad |z|<R.
\end{equation}
Choose $R>0$ small enough to avoid all other poles and zeros. A nonzero
norm twist may reduce $R$, but does not prevent this choice. After
decreasing $R$, Cauchy's estimate gives
$|C_{F,G,j}|\le C R^{-j}$.

For a Gevrey test $W$ constant near zero, set $f(u)=e^uW(u)$ and use
the Fourier convention of Lemma~\ref{lem:smoothed}. Its transform
satisfies $|\widehat f(t)|\le C e^{-a_0\sqrt{|t|}}$. Direct Fourier
inversion gives
\begin{equation}\label{eq:AlocalFourier}
 S_{F,G}(\lambda;W)=\frac1{2\pi}\int_\R\widehat f(t)
           A_{F,G}\!\left(1+\frac{1-it}{\lambda}\right)dt.
\end{equation}
Absolute convergence follows from Lemma~\ref{lem:bounds}; throughout
the line the Dirichlet series is bounded in absolute value by
$O(\lambda^2)$.

Expand the regular part in \eqref{eq:fullAjet} through degree $N-1$ on
$|t|\le\rho\lambda$, with fixed $0<\rho<R/4$. Its integrated remainder
is bounded by
\begin{equation}\label{eq:GevreyTaylorerror}
 \left(\frac{C_1}{\lambda}\right)^N(2N+2)!.
\end{equation}
Indeed the geometric Taylor remainder is
$O((|1-it|/(R\lambda))^N)$ on this interval, and the integral of
$(1+|t|)^N e^{-a_0\sqrt{|t|}}$ is at most $C_2^N(2N+2)!$.
The complementary integral of the original series is
$O(\lambda^{C_3}e^{-a_0\sqrt{\rho\lambda}})$.
When replacing the truncated polynomial integral by its whole-line
integral, use $\sum_{j<N}x^j\le N(1+x^N)$. The additional error is
bounded by the same exponential tail and $N$ times
\eqref{eq:GevreyTaylorerror}.

For every $j\ge0$, Fourier inversion gives the exact identity
\begin{equation}\label{eq:Taylorfunctional}
 \frac1{2\pi}\int_\R\widehat f(t)(1-it)^jdt
       =(-1)^jW^{(j)}(0).
\end{equation}
For example, the left side is $(1-d/du)^jf(0)$, and $f=e^uW$.
All terms with $j\ge1$ vanish for a plateau. The pole term gives
$r\lambda^2\int_0^\infty uW(u)\,du$, as in
Lemma~\ref{lem:smoothed}. Choose
$N=\lfloor\varepsilon\sqrt\lambda\rfloor$, with sufficiently small
fixed $\varepsilon>0$. Stirling's estimate in
\eqref{eq:GevreyTaylorerror} yields
\begin{equation}\label{eq:localGevreysum}
 S_{F,G}(\lambda;W)=r\lambda^2\int_0^\infty uW(u)\,du
          +C_{F,G}W(0)+O(e^{-a_1\sqrt\lambda}).
\end{equation}
Now apply Proposition~\ref{prop:reduction} with $W=\phi$.
The power saving in $T$ is smaller than
\eqref{eq:unconditionalpairerror}. This proves the theorem.
\end{proof}

The PNT-based proof remains useful for arbitrary smooth plateaus and
for its explicit treatment of norm twists. Theorem~\ref{thm:localGevrey}
is the cleanest fixed-representation version for the Gevrey tests.
Neither argument enlarges the Fourier support.

\section{\texorpdfstring{$L$}{L}-function probes and automorphic induction}
\label{sec:probes}

The second purpose of this paper is to use the mixed correlation as a
probe of the factorization \eqref{eq:induction}.  Throughout this section
that factorization is assumed to exist.  We write $F=L(s,\Pi)$ and
$P=L(s,\pi)$, where $\pi$ is a prescribed unitary cuspidal representation
of $\GL_{m_0}(\mathbb A_\Q)$.  Define
\[
 r_\Pi(\pi)=
 \begin{cases}
 r_\alpha,&\pi\simeq\sigma_\alpha\text{ for some }\alpha,\\
 0,&\text{otherwise}.
 \end{cases}
\]
The use of an inequivalent list of constituents makes this definition
unambiguous.  In particular, a positive answer identifies a constituent
of degree $m_0$ and determines its multiplicity.

Use the same $c,g,\phi$ in every pair correlation.  The zero multiset of
a product is the sum of the zero multisets of its factors.  Expanding
all ordered pairs therefore gives the exact identity
\begin{equation}\label{eq:productidentity}
 \frac12\{\mathcal R_{FP,FP}(T)-\mathcal R_{F,F}(T)
                         -\mathcal R_{P,P}(T)\}
 =\mathcal R_{F,P}(T).
\end{equation}
Coincident zeros and repeated factors cause no difficulty: their
multiplicities are retained on both sides.  The common normalization is
essential.  In particular, one must not rescale each product by its own
degree before taking this difference.

\subsection{Detecting a prescribed constituent}

For the constituent test, the Fourier amplitude can be chosen to remove
the smooth density altogether.  Let $\phi_*$ be a nonzero, nonnegative,
real, even Gevrey-order-two function supported in
$\{u: a_*\le |u|\le b_*\}$, with $0<a_*<b_*<b$ and $cb<1$.
Such functions are obtained by adding two reflected Gevrey bumps.
Then $\phi_*$ is a plateau of height zero near the origin and
\[
 J_{\phi_*}=\int_\R |u|\phi_*(u)\,du>0.
\]
We indicate the Fourier amplitude by a superscript when needed.

\begin{proposition}[The constituent probe]\label{prop:constituent-probe}
With the preceding choices, there is a constant $a>0$ such that
\begin{equation}\label{eq:probe-direct}
 \frac{\mathcal R^{\phi_*}_{FP,FP}(T)
       -\mathcal R^{\phi_*}_{F,F}(T)
       -\mathcal R^{\phi_*}_{P,P}(T)}{2J_{\phi_*}}
 =r_\Pi(\pi)+O\!\left(e^{-a\sqrt{\log T}}\right).
\end{equation}
Consequently the limit of the expression on the left is a nonnegative
integer.  It is positive precisely when $\pi$ is a constituent in
\eqref{eq:induction}, and then equals its multiplicity.
\end{proposition}

\begin{proof}
For a single cuspidal probe, \eqref{eq:degrees} gives
$r_{F,P}=r_\Pi(\pi)$ and $D_P=m_0$.  Apply
Theorem~\ref{thm:localGevrey} to $(F,P)$ with the zero plateau $\phi_*$.
Every smooth-density term in \eqref{eq:main} contains $\phi_*(0)$ and
therefore vanishes.  The finite-part term also contains $\phi_*(0)$.
Thus
\[
 \mathcal R^{\phi_*}_{F,P}(T)
 =r_\Pi(\pi)J_{\phi_*}
       +O\!\left(e^{-a\sqrt{\log T}}\right).
\]
Insert this in \eqref{eq:productidentity} and divide by $J_{\phi_*}$.
\end{proof}

This test requires no subtraction involving the conductor or the gamma
factors.  For a fixed probe its value is eventually within $1/2$ of the
integer $r_\Pi(\pi)$, so that rounding gives the multiplicity once the
error is sufficiently small.  This observation is an asymptotic
statement, not a uniform effective height bound for a family of probes.
A nonzero norm twist of a constituent is not counted unless it is
actually isomorphic to that constituent.  Fixed nonzero twist
parameters contribute only to the error in \eqref{eq:probe-direct}; the
constant in that error is not uniform when such a parameter tends to
zero.

\subsection{The arithmetic response after detecting the multiplicity}

A second Fourier amplitude, now with $\phi(0)\ne0$, retains the
secondary arithmetic response.  Set
\begin{equation}\label{eq:denssub}
 \mathcal D_{F,P}(T)=\mathcal R_{F,P}(T)
 -\frac{\phi(0)}{c^2L^2}
       \av{(D_FL+\eta_F)(m_0L+\eta_P)}.
\end{equation}
Theorems~\ref{thm:qualitative} and \ref{thm:localGevrey} give,
respectively,
\begin{equation}\label{eq:detect}
 \mathcal D_{F,P}(T)=r_\Pi(\pi)J_\phi
       +\frac{2\Rea C_{F,P}\,\phi(0)}{c^2L^2}+o(L^{-2}),
\end{equation}
and, for a Gevrey-order-two plateau, the same formula with remainder
$O(e^{-a\sqrt L})$.  If also $J_\phi\ne0$, this yields
\begin{align}
 r_\Pi(\pi)&=\lim_{T\to\infty}
                  \frac{\mathcal D_{F,P}(T)}{J_\phi},
                  \label{eq:multlimit}\\
 \Rea C_{F,P}&=\lim_{T\to\infty}\frac{c^2L^2}{2\phi(0)}
       \{\mathcal D_{F,P}(T)-r_\Pi(\pi)J_\phi\}.
                  \label{eq:Climit}
\end{align}
Alternatively, obtain $r_\Pi(\pi)$ first from
Proposition~\ref{prop:constituent-probe}; then the second limit only
requires $\phi(0)\ne0$.  A nonzero secondary response by itself does not
assert that $\pi$ is a constituent.  It is the leading coefficient in
\eqref{eq:probe-direct} that detects occurrence.

\subsection{Polarizing higher secondary terms}

The higher-level calculation supplies further probe responses.  Fix
$n$, use a common cutoff $h$ with $\int h^n>0$, and take a
permutation-invariant Fourier amplitude.  Write
$\mathcal C_n(F^{[n-j]},P^{[j]};T)$ for the raw mixed correlation with
$n-j$ entries from $F$ and $j$ entries from $P$, using the definition
\eqref{n:C-definition}.  The square brackets here indicate repeated
entries, not powers of a function.  For a nonnegative integer $z$,
expansion of the zero multiset of $FP^z$ gives
\begin{equation}\label{eq:probe-polarization}
 \mathcal C_n(FP^z;T)
 =\sum_{j=0}^n\binom nj z^j
      \mathcal C_n(F^{[n-j]},P^{[j]};T).
\end{equation}
There are $\binom nj$ choices of the probe positions, and each such
position has $z$ identical copies available.  This proves the identity
at finite height, before any asymptotic expansion.  Hence values at
$z=0,1,\ldots,n$ determine the mixed correlations by ordinary polynomial
interpolation.  The same identity applies coefficient by coefficient
to the expansion in Theorem~\ref{n:all-orders}.

There is an analogous polarization of the connected arithmetic
invariants.  The local expression $Y_{F,p}$ is additive under products,
and joint cumulants are multilinear.  Writing $\Gamma_n(H)=\Gamma(H,\ldots,H)$, we obtain, for $n\ge3$,
\begin{equation}\label{eq:probe-connected}
 \Gamma_n(FP^z)=\sum_{j=0}^n\binom nj z^j
       \Gamma(F^{[n-j]},P^{[j]}).
\end{equation}
Here the invariant on the right is the convergent mixed prime-power
cumulant defined in the higher-level calculation.  The identity also
holds for its frequency moments.  It makes explicit how repeated
product experiments separate responses that a single product combines.

For example, let
\[
 E_z=\Gamma_3(FP^z)-\Gamma_3(F)-z^3\Gamma_3(P).
\]
Then $E_z=3z\Gamma(F,F,P)+3z^2\Gamma(F,P,P)$.  The two choices $z=1,2$
give
\begin{align}
 \Gamma(F,F,P)&=\frac{4E_1-E_2}{6},\label{eq:probe-cubic-FFP}\\
 \Gamma(F,P,P)&=\frac{E_2-2E_1}{6}.\label{eq:probe-cubic-FPP}
\end{align}
In extracting $E_z$ from zero correlations, one first subtracts the
archimedean and pair contributions specified by the third-level
formula.  Thus these are arithmetic observables supplied by the proved
expansion, rather than an assertion that a raw third-order coefficient
consists only of $\Gamma_3$.

\subsection{A degree-one example}

Let $\chi$ be a primitive Dirichlet character of conductor $q$, viewed
as an automorphic representation of $\GL_1$. For the self coefficient
series,
\[
 A_{\chi,\chi}(s)
 =A_{\zeta,\zeta}(s)-\sum_{p\mid q}\frac{(\log p)^2}{p^s-1}.
\]
Therefore
\begin{equation}\label{eq:charconstant}
 C_{\chi,\chi}=C_{\zeta,\zeta}
                 -\sum_{p\mid q}\frac{(\log p)^2}{p-1}.
\end{equation}
With the Stieltjes-constant convention
\[
 \zeta(1+z)=\frac1z+\sum_{j\geq0}\frac{(-1)^j\gamma_j}{j!}z^j,
\]
formula \eqref{eq:Cbeta} gives
\begin{equation}\label{eq:zetaconstant}
 C_{\zeta,\zeta}=-\gamma_0^2-2\gamma_1
                   -\sum_p\frac{(\log p)^2}{(p-1)^2}.
\end{equation}
Indeed, $K=-\zeta'/\zeta$ has
$\beta_1=\gamma_0^2+2\gamma_1$, and
$\sum_{k\geq2}(k-1)p^{-k}=(p-1)^{-2}$.
This example shows arithmetic dependence beyond the degree and the
self multiplicity. It also shows that a single self secondary constant
is not a complete identifier: characters of the same conductor have
the same value in \eqref{eq:charconstant}.

\subsection{Scope of the probe and further development}

The statistic in Proposition~\ref{prop:constituent-probe} answers a
precise question: does a prescribed cuspidal representation occur in
the assumed factorization, and with what multiplicity?  Its degree is
already known from the choice of probe.  If one has a finite list known
to contain every constituent, applying the proposition to each member
recovers all multiplicities and hence the complete degree list.
More generally, no prior guarantee that the list is exhaustive is needed
if the inequivalent probes already detected account for the total degree:
\[
 \sum_{\pi\ \mathrm{detected}}r_\Pi(\pi)\deg(\pi)=D_F.
\]
Indeed every omitted constituent would contribute a positive amount to
the remaining degree. This is a stopping certificate once the
multiplicities have been determined, not a uniform finite-height
algorithm.
The present argument does not supply an exhaustive list, prove the
existence of automorphic induction, or recover every constituent from
finitely many arbitrary probes.

Self correlations and probe correlations have different roles.  The
leading self statistic records $D_F$ and $B_F$, while higher secondary
terms retain additional local arithmetic.  Nevertheless, as the
examples below show, even the full inverse-logarithmic self hierarchy
is not an injective invariant.  The mixed quantities
$\Rea C_{F,P}$, $\Gamma(F,F,P)$, and their higher analogues form a natural
collection of further responses to chosen probes.  Determining an
efficient or complete class of probes is a separate reconstruction
problem and a direction for further work.

Finally, the strict support restriction is unchanged.  Its decisive
consequence is that the relevant integer products have length
$X=T^{cb+o(1)}$ smaller than the height scale.  At the support boundary,
unequal products are no longer negligible by the present argument.
Improving the estimate for equal products does not by itself control
these new off-diagonal terms.

\section{The analytic reduction at an arbitrary fixed level}
\subsection{A balanced-packet reduction at every fixed level}
\label{n:analytic-reduction}

We give the analytic reduction with different height cutoffs. This
generality is useful when several coordinates are subsequently identified:
the cutoff attached to a block is the product of the original cutoffs.
Fix products $F_1,\ldots,F_n$ of unitary cuspidal $L$-functions over
$\mathbb Q$, and let $D_j,Q_j,b_j(p^k)$ denote their degrees,
conductors, and actual standard local coefficients. Write
\[
 c_j(N)=\Lambda(N)b_j(N),\qquad
 b_j^+(N)=\overline{b_j(N)},\qquad b_j^-(N)=b_j(N).
\]
The coefficients are zero unless $N$ is a prime power. The superscript
in $b_j^\pm$ specifies the sign of the Fourier-packet center, not the
sign preceding the arithmetic part of the explicit formula.

Choose $g_j\in C_c^\infty(\mathbb R)$ and put
\[
 h_j(t)=\int_{\mathbb R}g_j(v)e^{itv}\,dv,\qquad
 H(t)=\prod_{j=1}^n h_j(t),\qquad
 \kappa_{\boldsymbol h}=\int_{\mathbb R}H(t)\,dt\ne0.
\]
No positivity of $H$ is needed for the estimates. We use
\[
 \langle v\rangle_{\boldsymbol h}
  =\kappa_{\boldsymbol h}^{-1}\int_{\mathbb R}H(t)v(t)\,dt.
\]
Let $L=\log T$, $U=T/c$, and $\lambda=cL$. Suppose
\begin{equation}\label{n:strict-support}
 \Phi\in C_c^\infty(\mathbb R^n),\qquad
 \operatorname{supp}\Phi\subset
 \left\{u:\sum_{j=1}^n|u_j|<2b\right\},\qquad cb<1.
\end{equation}
Define the translation-invariant test and the raw correlation by
\begin{align}
 f(x)&=\int_{\mathbb R^n}\Phi(u)
     \delta(u_1+\cdots+u_n)e^{-2\pi i u\cdot x}\,du,
       \label{n:f-convention}\\
 \mathcal C_n(T;\Phi)&=
 \frac{2\pi}{\kappa_{\boldsymbol h}TL}
 \sum_{\gamma_1,\ldots,\gamma_n}
   \prod_{j=1}^n h_j(\gamma_j/U)
   f\left(\frac{\lambda\gamma_1}{2\pi},\ldots,
          \frac{\lambda\gamma_n}{2\pi}\right).
       \label{n:C-definition}
\end{align}
Zeros are counted with their multiplicities, and all tuples are retained.
The entire test functions make the definition meaningful without a
Riemann hypothesis.

For $I\subset\{1,\ldots,n\}$, let $\Phi_I$ be obtained from $\Phi$
by setting the coordinates outside $I$ equal to zero. Introduce the
balanced arithmetic functional
\begin{equation}\label{n:balanced-moment}
\begin{split}
 \mathfrak M_{I,\lambda}[W]
 ={}&\frac{(-1)^{|I|}}{\lambda^{|I|}}
 \sum_{\substack{N_i=p_i^{k_i}\ (i\in I),\ k_i\ge1\\
                 \epsilon_i\in\{+1,-1\}\\
                 \prod_{i\in I}N_i^{\epsilon_i}=1}}
 \prod_{i\in I}
       \frac{\Lambda(N_i)b_i^{\epsilon_i}(N_i)}{\sqrt{N_i}}
 \ W\left(\left(\frac{\epsilon_i\log N_i}{\lambda}\right)_{i\in I}
        \right).
\end{split}
\end{equation}
For a compactly supported $W$ this is a finite sum. Set
$\mathfrak M_{\varnothing,\lambda}[\Phi_\varnothing]=\Phi(0)$.
For $|I|=1$ the functional is zero. The minus sign in each arithmetic
term of the explicit formula is responsible for $(-1)^{|I|}$.

Finally put
\begin{equation}\label{n:singleton-density}
 \eta_j(t)=\log Q_j+D_j\log\frac{|t|}{2\pi c},\qquad
 a_j(t;\lambda)=\frac{D_jL+\eta_j(t)}{\lambda}
              =\frac{D_j}{c}+\frac{\eta_j(t)}{\lambda}.
\end{equation}

\begin{proposition}[Reduction to balanced prime powers]
\label{n:reduction-theorem}
For fixed $n$, fixed products and cutoffs, and
\eqref{n:strict-support}, there exists $\varepsilon>0$ such that
\begin{equation}\label{n:reduction-formula}
 \mathcal C_n(T;\Phi)=
 \sum_{I\subset\{1,\ldots,n\}}
 \left\langle\prod_{j\notin I}a_j(t;\lambda)
 \right\rangle_{\boldsymbol h}
 \mathfrak M_{I,\lambda}[\Phi_I]
 +O(T^{-\varepsilon}).
\end{equation}
One may use any fixed $0<\varepsilon<(1-cb)/2$. This analytic
reduction does not require Hypothesis H or a quantitative prime number
theorem. Those inputs enter when the balanced arithmetic functional is
expanded or evaluated.
\end{proposition}

\begin{proof}
We first record the coefficient estimates needed to sum all packet
errors. Rankin--Selberg theory, with the finite ramified correction for
the actual standard local roots, gives
\begin{equation}\label{n:c-one-two}
 \sum_{N\le y}|c_j(N)|\ll y,\qquad
 \sum_{N\le y}|c_j(N)|^2\ll y\log(2y).
\end{equation}
Indeed, the diagonal estimate
$\sum_{N\le y}\Lambda(N)|b_j(N)|^2\ll y$ and Cauchy--Schwarz
give the first assertion. The second follows by
$\Lambda(N)\le\log y$. The same reasoning applies to a finite product
by Cauchy--Schwarz over its finitely many constituents.

For any fixed collection of $r$ coordinates, let
\[
 A(M)=\sum_{N_1\cdots N_r=M}\prod_{i=1}^r|c_i(N_i)|.
\]
Repeated partial summation in \eqref{n:c-one-two} gives
\begin{equation}\label{n:convolution-lone}
 \sum_{M\le y}A(M)\ll y(\log(2y))^{r-1},\qquad
 \sum_{M\le y}\frac{A(M)}{\sqrt M}
       \ll \sqrt y(\log(2y))^{r-1}.
\end{equation}
The divisor bound and Cauchy--Schwarz inside each convolution give,
for every fixed $\nu>0$,
\begin{align}
 |A(M)|^2
 &\le d_r(M)\sum_{N_1\cdots N_r=M}\prod_i|c_i(N_i)|^2,
       \nonumber\\
 \sum_{M\le y}|A(M)|^2
 &\ll_\nu y^{1+\nu}(\log(2y))^{2r-1}.
       \label{n:convolution-ltwo}
\end{align}
These bounds also hold with conjugated coefficients, or with
coefficients replaced by smaller nonnegative quantities.

Let $\Omega_j$ be the exact gamma-factor density of $F_j$, and set
\[
 q_{j,U}(x)=\frac1{2\pi}\int_{\mathbb R}
               h_j(t)\Omega_j(Ut)e^{-ixt}\,dt.
\]
The gamma-factor estimates used in the pair case imply
\[
 \|h_j\Omega_j(U\cdot)\|_1\ll L,\quad
 \|(h_j\Omega_j(U\cdot))'\|_1\ll L,\quad
 \|(h_j\Omega_j(U\cdot))''\|_1\ll U+L.
\]
Consequently
\begin{equation}\label{n:q-estimates}
 |q_{j,U}(x)|\ll\frac{L}{1+|x|},\qquad
 |q_{j,U}(x)|\ll\frac{U+L}{x^2}\quad(x\ne0),\qquad
 \|q_{j,U}\|_1\ll L^2,\qquad
 \sup_x|xq_{j,U}(x)|\ll L.
\end{equation}
For the $L^1$ estimate, integrate the first bound up to $U$ and the
second beyond $U$. Thus no uniform estimate for arbitrarily high
derivatives of the scaled gamma factors is being assumed.

The explicit formula gives
\begin{equation}\label{n:EF-packets}
\begin{split}
 Z_j(u):={}&\sum_{\gamma_j}h_j(\gamma_j/U)e^{-i\lambda u\gamma_j}\\
 ={}&Uq_{j,U}(U\lambda u)+P_j(u)\\
 &-U\sum_{N\ge2}\frac{\Lambda(N)b_j(N)}{\sqrt N}
                      g_j(U\lambda u+U\log N)\\
 &-U\sum_{N\ge2}\frac{\Lambda(N)\overline{b_j(N)}}{\sqrt N}
                      g_j(U\lambda u-U\log N).
\end{split}
\end{equation}
The completed-function poles form a finite set, and
$|P_j(u)|\ll e^{\lambda|u|/2}$. Their fixed vertical shifts only
change the implied constant. By Fourier inversion,
$\mathcal C_n$ is $2\pi/(\kappa_{\boldsymbol h}U\lambda)$ times
\[
 \int\Phi(u)\delta(u_1+\cdots+u_n)\prod_jZ_j(u_j)\,du.
\]

First consider a term with arithmetic coordinates $I$, gamma
coordinates $J=I^c$, and no poles. Choose the arithmetic signs and
integers, and write
\[
 z_i=\epsilon_i\log N_i/\lambda\ (i\in I),\qquad
 z_j=0\ (j\in J),\qquad
 M=\prod_{\epsilon_i=+1}N_i,\quad
 N=\prod_{\epsilon_i=-1}N_i.
\]
The arithmetic packets have bounded Fourier width. On the support of
their product with $\Phi$, therefore,
\begin{equation}\label{n:product-length}
 MN=\prod_{i\in I}N_i
     \le \exp(2b\lambda+O(U^{-1}))=:X^2,\qquad
 X=T^{cb+o(1)}.
\end{equation}
This is a bound for the product of the arithmetic lengths; replacing
it by separate length bounds would lose the desired support range.

Make the change of variables
$y_j=U\lambda(u_j-z_j)$. The packet term has the factor
$U/\lambda^{n-1}$, and the delta constraint becomes
\[
 \sum_jy_j=-U\log(M/N).
\]
The prime coordinates contribute compactly supported $g_j(y_j)$;
the gamma coordinates contribute $q_{j,U}(y_j)$. Replacing
$\Phi(z+y/(U\lambda))$ by $\Phi(z)$ costs at most
$(U\lambda)^{-1}\|\nabla\Phi\|_\infty\sum_j|y_j|$.
This bound remains valid outside the support of either function, so
the replacement does not introduce a support exception.

For completeness, the weighted constrained convolution is bounded
uniformly in its constraint. If a gamma coordinate supplies the
factor $|y_j|$, eliminate that coordinate by the delta function and
use $\sup|y_jq_{j,U}(y_j)|\ll L$; every remaining gamma factor has
$L^1$ norm $O(L^2)$. If a prime coordinate supplies the weight, it is
bounded on its compact support, and one eliminates a gamma coordinate
using $\|q_{j,U}\|_\infty\ll L$. When there are no gamma coordinates,
all factors have compact support. Thus the constrained convolution
and its first absolute moment are bounded by a fixed power of $L$.
Equations \eqref{n:product-length} and \eqref{n:convolution-lone}
now show that the total normalized replacement error is
\begin{equation}\label{n:freeze-error}
 O\left(U^{-1}X L^{C_n}\right).
\end{equation}

After this replacement, convolution and Fourier inversion identify
the remaining constrained integral as
\begin{equation}\label{n:frozen-convolution}
 \frac1{2\pi}\int_{\mathbb R}H(t)
       \prod_{j\in J}\Omega_j(Ut)
       e^{iU\log(M/N)t}\,dt.
\end{equation}
At $M=N$ this is the desired common-height average. At $M\ne N$,
one integration by parts bounds it by
\begin{equation}\label{n:frozen-off-diagonal}
 \ll\frac{L^{|J|}}{U|\log(M/N)|}.
\end{equation}
The required derivative has $L^1$ norm $O(L^{|J|})$, as follows
by differentiating one factor at a time and integrating the bound
$U|\Omega'_j(Ut)|\ll U/(1+U|t|)$.
If $J$ is empty, the integral instead has compact Fourier support.
It vanishes exactly for all unequal $M,N$ in
\eqref{n:product-length}, once $T$ is large: when $M$ and $N$ are
comparable they are $O(X)$ and $U|\log(M/N)|\gg U/X\to\infty$;
when they are not comparable the same conclusion is immediate.

We spell out the sum of \eqref{n:frozen-off-diagonal}. In the region
where $M/N$ stays away from $1$, the reciprocal logarithm is bounded,
and \eqref{n:convolution-lone} gives $O(XL^{C_n})$ for the remaining
coefficient sum. In the region $M\asymp N\asymp V$, one has $V\ll X$
and
\[
 \frac1{\sqrt{MN}|\log(M/N)|}\ll\frac1{|M-N|}.
\]
For each nonzero shift $h$, Cauchy--Schwarz and
\eqref{n:convolution-ltwo} give
\[
 \sum_{M\asymp V}A(M)B(M+h)\ll_\nu V^{1+\nu}L^{C_n}.
\]
Summing $1/|h|$ and then summing the
dyadic ranges yields $O_\nu(X^{1+\nu}L^{C_n})$. If one sign set is
empty, its product is $1$, and every nonempty term lies in the first
region. The total normalized off-diagonal contribution is therefore
\begin{equation}\label{n:total-off-diagonal}
 O_\nu\left(U^{-1}X^{1+\nu}L^{C_n}\right).
\end{equation}

At $M=N$, Stirling's formula and direct integration near $t=0$ give
\begin{equation}\label{n:density-replacement}
 \int H(t)\left\{\prod_{j\in J}\Omega_j(Ut)
       -\prod_{j\in J}(D_jL+\eta_j(t))\right\}dt
       \ll U^{-1}L^{C_n}.
\end{equation}
On $|t|\le U^{-1}$ substitute $v=Ut$ and integrate the finite powers
of $|\log|v||$. On $|t|>U^{-1}$ use the error
$O((U|t|)^{-1})$ in each individual gamma factor and telescope the
product. The balanced absolute coefficient sum is
$O_\nu(X^\nu L^{C_n})$: group by $M=N$, divide by $M$, and apply
\eqref{n:convolution-ltwo} on dyadic intervals. Thus
\eqref{n:density-replacement} still contributes a power saving after
all balanced tuples are summed.

It remains to consider terms containing at least one pole coordinate.
Use $w_j=\lambda u_j$ and eliminate one pole coordinate by the delta
function. On the support of $\Phi$,
\[
 \exp\left(\frac12\sum_{j\text{ pole}}|w_j|\right)
 \le X\left(\prod_{i\text{ arithmetic}}N_i\right)^{-1/2}
       \exp\left(-\frac12\sum_{j\text{ gamma}}|w_j|\right)
       e^{O(U^{-1})}.
\]
Multiplication by the arithmetic coefficients therefore gives
$X\prod_i|c_i(N_i)|/N_i$. Each arithmetic packet integral cancels its
factor $U$; each gamma integral is bounded by
$\int U|q_{j,U}(Uw)|e^{-|w|/2}\,dw\le\|q_{j,U}\|_1\ll L^2$.
The remaining pole coordinates occupy intervals of length $O(L)$,
and partial summation in \eqref{n:c-one-two} gives
$\sum_{N\le X^2}|c_i(N)|/N\ll L$. Consequently the entire raw pole
contribution is $O(XL^{C_n})$, and its normalized value is
$O(U^{-1}XL^{C_n})$.

Finally, at a balanced tuple, \eqref{n:frozen-convolution} contributes
$\kappa_{\boldsymbol h}/(2\pi)$ times the common-height average.
The raw factor $U/\lambda^{n-1}$ and the normalization
$2\pi/(\kappa_{\boldsymbol h}U\lambda)$ leave $\lambda^{-n}$.
Assigning $\lambda^{-1}$ to each gamma singleton and each arithmetic
coordinate gives exactly \eqref{n:reduction-formula}. Choose $\nu>0$
small enough in \eqref{n:total-off-diagonal}, and absorb the fixed
powers of $L$ into a slightly smaller power saving. Since $cb<1$,
every error above is $O(T^{-\varepsilon})$ for any prescribed
$0<\varepsilon<(1-cb)/2$.
\end{proof}

\begin{remark}\label{n:cutoff-closure}
If coordinates in a block $B$ are identified, the new cutoff is
$h_B=\prod_{j\in B}h_j$. Its Fourier transform is the convolution of
the corresponding $g_j$, hence still belongs to $C_c^\infty$.
Furthermore the product of all block cutoffs is $H$, so
$\kappa_{\boldsymbol h}$ is unchanged. The same analytic reduction
therefore applies to every term of the finite partition sieve, with
the same normalization and a power-saving error.
\end{remark}

\section{Higher prime-power resonances and connected arithmetic measures}
% Analysis fragment: independent derivation of higher connected arithmetic terms.
\subsection{Higher local resonances and their convergence}
Let $F_j$ be any fixed finite collection of products of unitary cuspidal
$L$-functions over $\mathbb Q$, and let $b_j(p^k)$ denote their actual
standard logarithmic-derivative coefficients.  At a prime $p$ introduce
\[
 Y_{j,p}(z)=(\log p)\sum_{k\ge1}
 \frac{b_j(p^k)z^k+\overline{b_j(p^k)}z^{-k}}{p^{k/2}},
 \qquad |z|=1,
\]
where integration in $z$ is normalized Haar measure.  This series converges
absolutely at each fixed prime.  The sign of the arithmetic term in the
explicit formula is $-Y_{j,p}$.
For $B\subset[n]$, $|B|=r\ge3$, put
\[
 \Gamma_B=(-1)^r\sum_p\operatorname{cum}(Y_{j,p}:j\in B),
\]
where
\[
 \operatorname{cum}(Y_j:j\in B)=
 \sum_{\mathcal P\vdash B}(-1)^{|\mathcal P|-1}(|\mathcal P|-1)!
 \prod_{S\in\mathcal P}\int\prod_{j\in S}Y_j(z)\,dz.
\]
Here $\mathcal P\vdash B$ means that $\mathcal P$ is an unordered
set partition of $B$ into nonempty, pairwise disjoint blocks: their
union is $B$, and $|\mathcal P|$ is the number of blocks. For example,
$\{\{1,2\},\{3\}\}\vdash\{1,2,3\}$.
The same convention applies to $\mathcal Q\vdash B$ and to
partitions of $[n]$ below. No probabilistic assertion about the zeros
is being made: Haar integration is an algebraic device for selecting
zero total exponent.

Here is a sufficient convergence input supplied by the two assertions of
Jiang's Theorem~1.3 over $\mathbb Q$, in the fixed-data form \eqref{eq:Jianginput}.  For each
cuspidal constituent $\sigma$, some $\eta_0>0$ satisfies
\[
 \sum_{p\notin S}\sum_{k\ge2}
 \frac{|a_\sigma(p^k)|^2+\lambda_{\sigma\times\widetilde\sigma}(p^k)}
 {p^{k(1-\eta_0)}}<\infty. \tag{J+}
\]
The ordinary coefficients $\lambda$ provide a convenient control of high
powers of the prime coefficient $a_\sigma(p)$, although this is not the
only available route.
Indeed, at an unramified prime the formal Euler-factor identity is
\[
 \sum_{K\ge0}\lambda_{\sigma\times\widetilde\sigma}(p^K)t^K
 =\exp\left(\sum_{k\ge1}\frac{|a_\sigma(p^k)|^2}{k}t^k\right).
\]
All coefficients in the exponent are nonnegative.  The factor
$\exp(|a_\sigma(p)|^2t)$ alone contributes
$|a_\sigma(p)|^{2K}/K!$ to the coefficient of $t^K$.  Consequently
\[
 |a_\sigma(p)|^{2K}\le K!\lambda_{\sigma\times\widetilde\sigma}(p^K).
 \tag{H1K}
\]

Fix $r\ge3$.  A local resonance is a choice $\epsilon_j\in\{-1,1\}$,
$k_j\ge1$ for which $\sum_j\epsilon_j k_j=0$.  Write
$2K=\sum_jk_j$.  Then $K\ge2$ and $k_j\le K$ for every $j$.
There is a number $\eta>0$ such that
\[
 \sum_p\sum_{\substack{\epsilon_j=\pm1,\ k_j\ge1\\
                  \sum_j\epsilon_jk_j=0}}
 p^{\eta K}(\log p)^r
 \prod_{j=1}^r\frac{|b_j(p^{k_j})|}{p^{k_j/2}}<\infty. \tag{Rr}
\]
To prove this, expand each $b_j$ into its finitely many constituents.
Let $\varepsilon_*>0$ be a common local gap, so
$|a_\sigma(p^k)|\le d_\sigma p^{k(1/2-\varepsilon_*)}$.
Choose $0<\eta'<\min(\eta_0,2\varepsilon_*)$ and first fix $K$.
Set
\[
 x_j=|a_{\sigma_j}(p^{k_j})|p^{-k_j(1-\eta')/2}.
\]
Weighted arithmetic--geometric mean gives
\[
 p^{\eta'K}\prod_j\frac{|a_{\sigma_j}(p^{k_j})|}{p^{k_j/2}}
 =\prod_jx_j
 \le\sum_j\frac{k_j}{2K}x_j^{2K/k_j}.
\]
For $k_j\ge2$, the exponent $2K/k_j$ is at least two and the local bound
makes $x_j$ bounded; hence its summand is bounded, for fixed $K$, by a
constant times $|a_{\sigma_j}(p^{k_j})|^2p^{-k_j(1-\eta')}$.
For $k_j=1$, inequality (H1K) bounds it by a constant times
$\lambda_{\sigma_j\times\widetilde\sigma_j}(p^K)p^{-K(1-\eta')}$.
Thus (J+) proves convergence for every fixed $K$.
Choose $0<\eta<\eta'$; the spare factor $p^{-(\eta'-\eta)K}$ absorbs
$(\log p)^r$.  This proves (Rr) for bounded $K$.
For large $K$, use the local bound directly.  There are $O_r(K^{r-1})$
compositions and finitely many signs, and their total is bounded by
\[
 C_r\sum_p(\log p)^r\sum_{K\ge K_0}
 K^{r-1}p^{-(2\varepsilon_*-\eta)K},
\]
which converges for sufficiently large $K_0$.
The finitely many ramified primes are treated by their strict local
bounds, decreasing $\eta$ if necessary.
This proves (Rr), and therefore absolute convergence of $\Gamma_B$,
all of its logarithmically weighted analogues, and a positive power tail
when the truncation variable is $p^K$.

For clarity, the coefficient of a particular resonant monomial in the
cumulant is
\[
 c(\epsilon,k)=
 \sum_{\substack{\mathcal P\vdash B\\
          \sum_{j\in S}\epsilon_j k_j=0\ (S\in\mathcal P)}}
       (-1)^{|\mathcal P|-1}(|\mathcal P|-1)!.
\]
It is bounded in terms of $r$ alone.  Thus the above proof controls all
cumulant subtractions as well as the moments themselves.

\begin{remark}
The $\lambda$-coefficient input in (J+) is a convenient sufficient input,
not a necessary new hypothesis.  The Rankin--Selberg estimate already used
in the pair calculation gives
$\sum_{p\le x}(\log p)|a_\sigma(p)|^2\ll x$.
Combining this with $|a_\sigma(p)|\ll p^\theta$, $\theta<1/2$, yields,
for $K\ge2$,
\[
 \frac{|a_\sigma(p)|^{2K}}{p^{K(1-\eta)}}
 \ll_K\frac{|a_\sigma(p)|^2}
 {p^{1+(K-1)(1-2\theta)-K\eta}}.
\]
If $0<\eta<(1-2\theta)/2$, the exponent on the right is strictly
larger than one for every $K\ge2$.  Partial summation therefore proves
convergence with any fixed logarithmic weight.  This supplies an
alternative to (H1K) in the preceding proof.  Thus the joint weighted
prime-power convergence furnished by Jiang's logarithmic-derivative
assertion, together with the usual Rankin--Selberg prime mean, already
suffices for the higher-resonance argument.
\end{remark}

\subsection{The first new connected constants}
For three labels, no proper partition has every block balanced, so
\[
 \Gamma_{123}=-2\operatorname{Re}\sum_p(\log p)^3
 \sum_{i=1}^3\sum_{a,b\ge1}
 \frac{b_j(p^a)b_k(p^b)\overline{b_i(p^{a+b})}}{p^{a+b}},
 \quad\{j,k\}=\{1,2,3\}\setminus\{i\}.
\]
In the self case this becomes
\[
 \Gamma_3(F)=-6\operatorname{Re}\sum_p(\log p)^3
       \sum_{a,b\ge1}\frac{b_F(p^a)b_F(p^b)
           \overline{b_F(p^{a+b})}}{p^{a+b}}.
\]
These terms come from $p^ap^b=p^{a+b}$, including $p\,p=p^2$.
Their first contribution to a normalized correlation has order
$\lambda^{-3}$, not $\lambda^{-2}$.

At level four the subtraction of the three pairings is essential:
\[
 \operatorname{cum}(Y_1,Y_2,Y_3,Y_4)
 =\int Y_1Y_2Y_3Y_4
  -\sum_{\{ij\}\{kl\}}\left(\int Y_iY_j\right)
                         \left(\int Y_kY_l\right).
\]
Write $b_k=b_F(p^k)$.  In the self case,
\begin{align*}
 \Gamma_4(F)=\sum_p(\log p)^4\bigg\{
 &8\operatorname{Re}\sum_{a,b,c\ge1}
       \frac{b_ab_bb_c\overline{b_{a+b+c}}}{p^{a+b+c}}\\
 &+6\sum_{\substack{a,b,c,d\ge1\\a+b=c+d}}
       \frac{b_ab_b\overline{b_c b_d}}{p^{a+b}}
 -12\left(\sum_{k\ge1}\frac{|b_k|^2}{p^k}\right)^2\bigg\}.
\end{align*}
In particular, the part in which every exponent is one is
$-6\sum_p(\log p)^4|b_F(p)|^4/p^2$.
This term must not be discarded: its coefficients involve the fourth
power of the prime coefficient.  Either the ordinary Rankin--Selberg
coefficient bound in (J+) or the prime-mean argument below supplies the
required estimate.

For a primitive Dirichlet character $\chi$ of conductor $q$, these formulas
have the particularly simple evaluations
\[
 \Gamma_3(\chi)=-6\sum_{p\nmid q}\frac{(\log p)^3}{(p-1)^2},
 \qquad
 \Gamma_4(\chi)=\sum_{p\nmid q}
       \frac{(26-6p)(\log p)^4}{(p-1)^3}.
\]
For the second identity, use
\[
 \sum_{k\ge1}p^{-k}=\frac1{p-1},\quad
 \sum_{a,b,c\ge1}p^{-a-b-c}=\frac1{(p-1)^3},\quad
 \sum_{a+b=c+d}p^{-a-b}=\frac{p+1}{(p-1)^3}.
\]
These are explicit arithmetic invariants beyond the leading degree and
multiplicity data.  They are not complete identifiers: characters of the
same conductor have the same displayed self constants.

\section{The full inverse-logarithmic expansion}
% Independent algebraic fragment. Analytic height reduction is in the audit fragment.
\subsection{The frequency-resolved cumulant expansion}
Use the notation of Proposition~\ref{n:reduction-theorem} and put
$a_j=D_j/c$. For the common-cutoff specialization write
$\langle v\rangle_n=(\int h^n)^{-1}\int h^n v$. One may ensure
$\int h^n>0$ simultaneously for every $n$ by taking $h=h_0^2$, where
$h_0$ is real, even, and has a smooth compactly supported Fourier transform.
The general formulas use $\langle\cdot\rangle_{\boldsymbol h}$.

For $B\subset[n]$, define a local complex frequency measure $m_{B,p}$ by
\[
 m_{B,p}=(-1)^{|B|}(\log p)^{|B|}
 \sum_{\substack{\epsilon_j=\pm1,\ k_j\ge1\\\sum_{j\in B}\epsilon_jk_j=0}}
 \prod_{j\in B}\frac{b_j^{[\epsilon_j]}(p^{k_j})}{p^{k_j/2}}
 \delta_{(-\epsilon_jk_j\log p)_{j\in B}},
 \quad b_j^{[+]}=b_j,\quad b_j^{[-]}=\overline{b_j}.
\]
In particular $m_{\{j\},p}=0$. Its cumulant measure is
\[
 \nu_{B,p}=\sum_{\mathcal Q\vdash B}
 (-1)^{|\mathcal Q|-1}(|\mathcal Q|-1)!
 \bigotimes_{A\in\mathcal Q}m_{A,p},\qquad \nu_B=\sum_p\nu_{B,p}.
\]
The cumulant subtraction is essential: two pair blocks at the same prime
must not be treated as independent. The exact prime-power resonance
moment is the sum over set partitions into blocks of size at least two
of tensor products of the $\nu_B$. This identity follows first with a
finite prime cutoff, and then on compactly supported frequency tests.

For $|B|\ge3$ put
\[
 M_B(\alpha)=\int y^\alpha\,d\nu_B(y),\qquad
 \Gamma_B=M_B(0).
\]
The positive power moment proved above supplies absolute convergence
of these moments and every further moment required for a remainder estimate.
The normalized block operator has expansion
\[
 \mathcal H_B(\lambda)\Phi
 =\sum_{\alpha}\lambda^{-|B|-|\alpha|}
 \frac{M_B(\alpha)}{\alpha!}
 (\partial_B^\alpha\Phi)|_{\xi_B=0}.
\]
This is an asymptotic operator expansion on compact sets of test functions.

For a pair $e=\{i,j\}$ choose $i<j$ and $v_e=e_i-e_j$. Suppose
\[
 A_{ij}(1+z)=\frac{r_{ij}}{z^2}+\sum_{m\ge0}C_{ij,m}z^m,
 \qquad C_{ij,0}=C_{ij}.
\]
Then the pair operator is
\[
 \mathcal P_{ij}(\lambda)\Phi
 =r_{ij}\int_{\mathbb R}|u|\Phi(uv_e)\,du
 +\sum_{m\ge0}\lambda^{-m-2}
 \{C_{ij,m}+(-1)^m\overline{C_{ij,m}}\}
 (\partial_{v_e}^m\Phi)|_{\xi_i=\xi_j=0},
\]
with coordinates outside the pair left free in both terms. This sign
convention places $b_i\overline{b_j}$ on the negative $v_e$ ray, as in
the explicit formula. Indeed, for a one-sided test $W$ the corresponding
smoothed sum is
\[
 r_{ij}\lambda^2\int_0^\infty uW(u)\,du
 +\sum_{m\ge0}(-1)^m C_{ij,m}\lambda^{-m}W^{(m)}(0).
\]
For permutation-symmetric $\Phi$, all odd pair derivatives vanish.

Each singleton $j$ contributes the operator
\[
 \mathcal A_j(t,\lambda)\Phi
 =\left(a_j+\frac{\eta_j(t)}\lambda\right)\Phi|_{\xi_j=0}.
\]
\begin{theorem}[Expansion at every fixed level and order]
\label{n:all-orders}
Let the fixed products, cutoffs, and Fourier amplitude satisfy
Proposition~\ref{n:reduction-theorem}. For a block $B$ of a partition,
let $\mathcal T_B$ denote $\mathcal A_j$ if $B=\{j\}$,
$\mathcal P_{ij}$ if $B=\{i,j\}$, and $\mathcal H_B$ if $|B|\ge3$.
For every fixed nonnegative integer $M$,
\begin{equation}\label{n:master-expansion}
 \mathcal C_n(T;\Phi)=
 \left[
 \left\langle
 \sum_{\mathcal P\vdash[n]}
     \left(\bigotimes_{B\in\mathcal P}\mathcal T_B(t,\lambda)\right)
                  \Phi
 \right\rangle_{\boldsymbol h}
 \right]_{\le M}
 +O(\lambda^{-M-1})+O(T^{-\varepsilon}).
\end{equation}
The sum $\mathcal P\vdash[n]$ runs over all unordered set
partitions of $[n]$, with the same convention as above.
The brackets retain terms of total inverse-$\lambda$ degree at most
$M$. In particular, this is a finite explicit prescription for every
coefficient. For a larger block one retains $|B|+|\alpha|\le M$;
for a pair, one retains $m+2\le M$, as well as its leading integral.
All displayed operators act only on their own coordinates, leaving the
other coordinates free until their operators are applied.
\end{theorem}

\begin{proof}
At a finite prime cutoff, introduce independent circle variables $z_p$
only as an algebraic device. The Haar integral of a monomial is zero
unless its exponent at each prime is zero. Expanding the product of
the sums of local arithmetic expressions therefore produces exactly
the balanced functional in \eqref{n:balanced-moment}.
The elementary moment--cumulant identity states that a joint moment
is the sum, over set partitions, of products of joint cumulants.
One proof is to expand the formal power series
\[
 \log\!\left(\int\exp\left(\sum_jt_jX_j\right)\right):
 \quad [t_1\cdots t_s]\log(1+V)
 =\sum_{\mathcal Q\vdash[s]}
       (-1)^{|\mathcal Q|-1}(|\mathcal Q|-1)!
       \prod_{B\in\mathcal Q}\int\prod_{j\in B}X_j.
\]
For variables attached to different primes the generating functions
multiply and their logarithms add. Thus the cumulant of a sum of the
local variables is the sum of their cumulants. The same finite algebra,
with the coordinate frequencies retained, proves the identity for the
measures $\nu_B$. There is no assertion here of independence of zeros.

For pairs, Proposition~\ref{prop:sufficient} and the local
Rankin--Selberg factorization give \eqref{eq:fullAjet} for every
$(F_i,F_j)$. Apply Fourier localization as in
\eqref{eq:AlocalFourier}, now with a fixed Taylor degree and a smooth
test which need not be a plateau. The tail of its Fourier transform
is smaller than any prescribed inverse power of $\lambda$, and
\eqref{eq:Taylorfunctional} gives the stated pair operator, including
its signs. The remainder after the term $m=M-2$ is
$O(\lambda^{-M-1})$, after the normalization $\lambda^{-2}$.

For $|B|\ge3$, the positive power moment proved in the preceding
section implies
\[
 \int e^{a_B\sum_{j\in B}|y_j|}\,d|\nu_B|(y)<\infty
\]
for some $a_B>0$. Taylor's formula with remainder gives
\[
 \lambda^{-|B|}\int W(y/\lambda)\,d\nu_B(y)
 =\sum_{|\alpha|\le M-|B|}
    \frac{M_B(\alpha)}{\alpha!\lambda^{|B|+|\alpha|}}
        \partial^\alpha W(0)+O(\lambda^{-M-1}).
\]
Every further moment required for this remainder is absolutely finite.
These estimates are uniform when the finite collection of test
seminorms in the proof is bounded.

The normalized pair measures have uniformly bounded total variation
on fixed compact frequency sets, by Lemma~\ref{lem:bounds}.
The larger-block measures have total variation $O(\lambda^{-|B|})$.
Consequently their tensor products can be expanded successively with
the same error bounds. Integration over a pair variable, restriction
to zero coordinates, and differentiation all preserve smoothness and
compact support in the remaining coordinates. There are only finitely
many partitions for fixed $n$. The singleton factors are exact
polynomials of degree one in $\lambda^{-1}$, and all products of the
functions $\eta_j$ are integrable against $H(t)$.
Combining this arithmetic expansion with
\eqref{n:reduction-formula} proves \eqref{n:master-expansion}.
\end{proof}

All frequency measures in this proof are supported on the appropriate
zero-sum hyperplanes. The derivative tensors of the larger blocks are
therefore tangent to those hyperplanes: although multi-index notation
uses an extension of $\Phi$ to the ambient space, the contracted result
does not depend on that extension. No additional copy of
$\delta(\sum u_j)$ is inserted when these measures are paired with
$\Phi$.

\subsection{The first two secondary orders}
Let $\mathcal M$ run over matchings of $[n]$, and let $U(\mathcal M)$
be its unmatched vertices. Put
\[
 I_{\mathcal M}(\Phi)=\int_{\mathbb R^{|\mathcal M|}}
 \prod_{e\in\mathcal M}|u_e|\,
 \Phi\left(\sum_{e\in\mathcal M}u_ev_e\right)\,du,
\]
and
\[
 A_U^{(s)}=\sum_{\substack{J\subset U\\|J|=s}}
 \left\langle\prod_{j\in J}\eta_j\right\rangle_{\boldsymbol h}
 \prod_{j\in U\setminus J}a_j.
\]
Then
\begin{align*}
 \mathcal C_n(T)={}&
 \sum_{\mathcal M}\left(\prod_{e\in\mathcal M}r_e\right)
 \left(A_U^{(0)}+\frac{A_U^{(1)}}\lambda+
                    \frac{A_U^{(2)}}{\lambda^2}\right)I_{\mathcal M}(\Phi)\\
 &+\frac1{\lambda^2}\sum_{\mathcal M}\sum_{e\in\mathcal M}
 2\operatorname{Re}C_e\left(\prod_{f\in\mathcal M\setminus\{e\}}r_f\right)
 A_U^{(0)}I_{\mathcal M\setminus\{e\}}(\Phi)
 +O(\lambda^{-3})+O(T^{-\varepsilon}).
\end{align*}
Genuine three-coordinate resonances first occur at order $\lambda^{-3}$;
genuine four-coordinate connected terms first occur at $\lambda^{-4}$.

In the self case write $a=D_F/c$, $B=B_F$, $C=C_{F,F}$,
$\mu_j=\langle\eta_F^j\rangle_n$, $q=n-2r$, and
$N_{n,r}=n!/(2^rr!(n-2r)!)$. For symmetric $\Phi$ this simplifies to
\[
 \mathcal C_n(T)=\sum_{r\le n/2}N_{n,r}B^rI_{n,r}(\Phi)
 \left\{a^q+\frac{qa^{q-1}\mu_1}\lambda
 +\frac{\binom q2a^{q-2}(\mu_2+2\operatorname{Re}C)}{\lambda^2}\right\}
 +O(\lambda^{-3})+O(T^{-\varepsilon}),
\]
where terms with zero combinatorial coefficient are omitted, so negative
powers of $a$ need not be interpreted.

\subsection{An explicit third-level formula}
For three possibly different products and a permutation-symmetric test,
let $\mu_S=\langle\prod_{j\in S}\eta_j\rangle_{\boldsymbol h}$. If $\{i,j,k\}=[3]$,
\begin{align*}
 \mathcal C_3(T)={}&\Phi(0)\bigg[
 a_1a_2a_3+\frac1\lambda\sum_i\mu_i\prod_{j\ne i}a_j
 +\frac1{\lambda^2}\sum_{i<j}(\mu_{ij}+2\operatorname{Re}C_{ij})a_k\\
 &\hspace{32mm}+\frac1{\lambda^3}
 \left\{\mu_{123}+\sum_{i<j}2\operatorname{Re}C_{ij}\mu_k
                   +\Gamma_{123}\right\}\bigg]\\
 &+\sum_{i<j}r_{ij}\left(a_k+\frac{\mu_k}\lambda\right)I_{ij}(\Phi)
 +O(\lambda^{-4})+O(T^{-\varepsilon}).
\end{align*}
Thus the self formula is
\begin{align*}
 \mathcal C_3(T)={}&\Phi(0)\left[
 a^3+\frac{3a^2\mu_1}\lambda
 +\frac{3a(\mu_2+2\operatorname{Re}C)}{\lambda^2}
 +\frac{\mu_3+6\operatorname{Re}C\mu_1+\Gamma_3(F)}{\lambda^3}\right]\\
 &+3B\left(a+\frac{\mu_1}\lambda\right)I_{3,1}(\Phi)
 +O(\lambda^{-4})+O(T^{-\varepsilon}).
\end{align*}
The first connected arithmetic invariant is
\[
 \Gamma_3(F)=-6\operatorname{Re}\sum_p(\log p)^3\sum_{a,b\ge1}
 \frac{b_F(p^a)b_F(p^b)\overline{b_F(p^{a+b})}}{p^{a+b}}.
\]
It comes from the exact relations $p^ap^b=p^{a+b}$, and therefore cannot
be obtained by multiplying pair-correlation constants alone.

At fourth order one must also retain the second derivative of the pair
regular part ($C_{ij,2}$), the connected fourth constant $\Gamma_4$,
and products of two finite pair constants. A pair-only subexponential
remainder cannot be transferred to higher levels after truncation at
order two or three: these explicit further inverse powers remain.

\section{A fourth-level formula and the distinct-index sieve}

\subsection{The self fourth-level formula through order four}

For this explicit display, take $F_j=F$, a common cutoff $h$, and a
real permutation-symmetric Fourier amplitude. Put
\[
 a=D_F/c,\qquad a(t;\lambda)=a+\eta_F(t)/\lambda,\qquad
 B=B_F,\qquad c_0=2C_{F,F},\qquad c_2=2C_{F,F,2}.
\]
The subscript $2$ denotes the coefficient of $z^2$ in
\eqref{eq:fullAjet}, not the second derivative itself. The coefficients
of this self series are real. Thus every odd pair derivative vanishes
even before symmetry of the test is used.

Write $\partial_{12}=\partial_{u_1}-\partial_{u_2}$, and define
\begin{align*}
 I_1&=\int_\R |v|\Phi(v,-v,0,0)\,dv,\\
 I_2&=\int_{\R^2}|vw|\Phi(v,-v,w,-w)\,dv\,dw,\\
 J_1&=\int_\R |v|\,
             (\partial_{12}^2\Phi)(0,0,v,-v)\,dv.
\end{align*}
The average in the next formula is against $h^4/\int h^4$.
Theorem~\ref{n:all-orders} gives
\begin{equation}\label{n:fourth-explicit}
\begin{split}
 \mathcal C_4(T)={}&\Phi(0)\left[
   \langle a(t;\lambda)^4\rangle_4
   +\frac{6c_0}{\lambda^2}\langle a(t;\lambda)^2\rangle_4
   +\frac{4\Gamma_3(F)}{\lambda^3}\langle a(t;\lambda)\rangle_4
   +\frac{3c_0^2+\Gamma_4(F)}{\lambda^4}\right]\\
 &+6B I_1\langle a(t;\lambda)^2\rangle_4
     +3B^2I_2+\frac{6Bc_0}{\lambda^2}I_1\\
 &+\frac{6c_2}{\lambda^4}
       \left\{a^2\partial_{12}^2\Phi(0)+BJ_1\right\}
     +O(\lambda^{-5})+O(T^{-\varepsilon}).
\end{split}
\end{equation}
This formula includes all terms through degree four. To check the
integer factors, list the five partition types of four labels:
\[
\begin{array}{c|ccccc}
 \text{type}&1+1+1+1&2+1+1&2+2&3+1&4\\\hline
 \text{number}&1&6&3&4&1.
\end{array}
\]
For type $2+1+1$, the pair supplies its leading integral, its constant
$c_0/\lambda^2$, or its second derivative
$c_2\partial_{12}^2/\lambda^4$. For type $2+2$, choosing the finite
constant in one pair gives $3\cdot2=6$ terms; choosing it in both
gives $3c_0^2\Phi(0)/\lambda^4$. Choosing the second derivative in one
pair gives the term $6Bc_2J_1/\lambda^4$.
The four triples give $4\Gamma_3(F)\langle a(t;\lambda)\rangle_4/
\lambda^3$. Their linear Taylor terms vanish: at coincident triple
coordinates, permutation symmetry makes the three first derivatives
equal, while their frequency sum is zero. The full block gives
$\Gamma_4(F)/\lambda^4$.

Thus $\Gamma_4$ is not the sole fourth-order contribution. In
particular, deleting same-prime collisions from the pair product without
the cumulant subtraction, or dropping the pair Taylor coefficient
$C_{F,F,2}$, would give an incorrect fourth-level expansion.

\subsection{Why the pair subexponential remainder does not transfer directly}

At level two a plateau kills every nonconstant Taylor term of the pair
regular part. At higher levels the constants $\Gamma_3,\Gamma_4,\ldots$
remain even when the Fourier amplitude is constant near the origin.
Moreover, a plateau near the origin alone need not kill pair derivatives
when other pair variables are nonzero. Therefore an expansion retaining
only the leading term and the coefficient of $(\log T)^{-2}$ cannot
have a subexponential remainder at every level. The omitted cubic term
already gives a counterexample. Formula~\eqref{n:master-expansion}
specifies the terms that must first be subtracted.

For each fixed desired order the remainder is controlled as stated in
Theorem~\ref{n:all-orders}. The present result does not assert convergence
of the resulting formal inverse-logarithmic series, nor does it use an
unproved interchange between the limits $n\to\infty$ and $T\to\infty$.
Both $n$ and the truncation order remain fixed.

\subsection{Removing repeated indices}

For a self correlation, the distinction between raw tuples and tuples
of distinct indices is exactly combinatorial. Let $\mathcal P_n$ be
the set of partitions of $[n]$. For $\mathcal F\in\mathcal P_n$, order
its blocks by their least elements, write $r=|\mathcal F|$, and put
\[
 \mu(\mathcal F)=\prod_{S\in\mathcal F}
                    (-1)^{|S|-1}(|S|-1)!,\qquad
 h_S(t)=\prod_{j\in S}h_j(t),\qquad
 f_{\mathcal F}(x)=f(y),\quad y_j=x_S\quad(j\in S).
\]
The exact sieve, in the normalization used here, is
\begin{equation}\label{n:exact-sieve}
 \mathcal C_n^{\rm dist}(T;f,\boldsymbol h)
   =\sum_{\mathcal F\in\mathcal P_n}\mu(\mathcal F)
      \mathcal C_r(T;f_{\mathcal F},(h_S)_{S\in\mathcal F}).
\end{equation}
An index includes its occurrence in the zero multiset. Thus different
copies of a repeated factor are distinct indices even when their
ordinates coincide. Formula~\eqref{n:exact-sieve} removes equality of
indices, not equality of their numerical ordinates.

For clarity, the Fourier amplitude to use in the lower-level term is
the smooth pushforward
\begin{equation}\label{n:pushforward}
 \Phi_{\mathcal F}(v)=\int_{\R^n}\Phi(u)
       \prod_{S\in\mathcal F}
             \delta\!\left(v_S-\sum_{j\in S}u_j\right)du.
\end{equation}
Indeed, direct substitution gives
\[
 f_{\mathcal F}(x)=\int_{\R^r}\Phi_{\mathcal F}(v)
       \delta\!\left(\sum_Sv_S\right)e^{-2\pi i v\cdot x}\,dv.
\]
Integration along the linear fibers in \eqref{n:pushforward} preserves
smoothness and compact support, and
$\sum_S|v_S|\le\sum_j|u_j|$. Thus it preserves the strict support
range. Also $\prod_Sh_S=\prod_jh_j$, so the normalizing $\kappa$ is
unchanged. The finite sieve therefore preserves every fixed-order
asymptotic expansion and its stated remainder.

To prove the coefficient, fix a tuple with equality partition
$\mathcal E$. Its weight on the right of \eqref{n:exact-sieve} is
$\sum_{\mathcal F\le\mathcal E}\mu(\mathcal F)$, which factors over
the blocks of $\mathcal E$. On a block of size $s$, the factor is the
sum of the signs of all permutations of $s$ letters, after grouping a
permutation by its disjoint cycles. This is one if $s=1$ and zero if
$s\ge2$. Hence the weight is one exactly for distinct indices.
This is the Rudnick--Sarnak sieve used in \cite{RS,LY02}; see also
\cite{WY} for its detailed determinant interpretation.

For $n=3$ the identity is
\[
 \mathcal C_3^{\rm dist}
 =\mathcal C_3-\mathcal C_2^{12|3}
                   -\mathcal C_2^{13|2}-\mathcal C_2^{23|1}
                   +2\mathcal C_1^{123}.
\]
The contracted pair and one-level terms contain no cubic prime-power
cumulant. Thus the new invariant $\Gamma_3$ survives the sieve; its
coefficient is still $\Phi(0)/\lambda^3$ in the raw three-variable
normalization. The lower terms must nevertheless be included when
reporting the full distinct-index secondary expansion.

\section{Arithmetic information beyond pair correlations}\label{sec:information}
\subsection{A strict gain over all pair arithmetic data}
The following explicit example shows that higher secondary coefficients
can distinguish products whose entire self-pair arithmetic data agree.
All factors have degree one, so this example does not depend on Jiang's
theorem or on any unproved higher-degree estimate.

Let $\chi$ be the primitive character modulo $13$ defined by
$\chi(2)=\zeta_{12}=e^{2\pi i/12}$; the residue class of $2$ generates
$(\mathbb Z/13\mathbb Z)^\times$.  Put
\[
 A=\{2,3,6,8\},\qquad B=\{1,2,4,8\}\subset\mathbb Z/12\mathbb Z,
 \qquad
 F_A=\prod_{a\in A}L(s,\chi^a),\quad
 F_B=\prod_{b\in B}L(s,\chi^b).
\]
Every exponent is nonzero modulo $12$, so all eight factors here are
primitive characters of conductor $13$.  Each product has exactly one
odd character and three even characters.  Consequently the two products
have the same degree $4$, constituent multiplicity invariant $4$,
conductor $13^4$, and complete archimedean gamma factors.

For $X=A,B$ and $d\in\mathbb Z/12\mathbb Z$, let
$N_X(d)=\#\{(x,y)\in X^2:x-y=d\}$.  Direct counting gives
\[
\begin{array}{c|rrrrrrrrrrrr}
 d&0&1&2&3&4&5&6&7&8&9&10&11\\\hline
 N_A(d)&4&1&1&1&1&1&2&1&1&1&1&1\\
 N_B(d)&4&1&1&1&1&1&2&1&1&1&1&1
\end{array}
\]
Thus, whenever $z^{12}=1$,
\[
 \left|\sum_{a\in A}z^a\right|^2
 =\sum_dN_A(d)z^d
 =\sum_dN_B(d)z^d
 =\left|\sum_{b\in B}z^b\right|^2.
\]
At $p\ne13$ take $z=\chi(p)^k$.  At $p=13$ all standard local
coefficients of both products vanish.  It follows that
\[
 |b_{F_A}(p^k)|^2=|b_{F_B}(p^k)|^2
 \quad\hbox{for every prime $p$ and every $k\ge1$}.
\]
Their self-pair coefficient Dirichlet series are therefore identical,
as are all its Laurent coefficients at $s=1$.  This conclusion is
stronger than equality of the first pair secondary constant.

The corresponding triple coefficient products differ.  At $p=2$,
using $\zeta_{12}=(\sqrt3+i)/2$, direct calculation gives
\[
\begin{array}{c|cc}
 &F_A&F_B\\\hline
 b_F(2)&-1+i&\dfrac{\sqrt3-1}{2}
                       +i\dfrac{1+\sqrt3}{2}\\[5pt]
 b_F(4)&-1+i\sqrt3&-1+i\sqrt3\\[3pt]
 b_F(2)^2\overline{b_F(4)}&-2\sqrt3+2i&2\sqrt3+2i.
\end{array}
\]
In particular, the real parts differ by $-4\sqrt3$.

For a proof that this difference is detectable by the secondary
hierarchy, let $\nu_{3,F}$ be the complex triple-resonance measure before
rescaling the frequency variables by $\lambda$.  It has atoms at
\[
 (a\log p,b\log p,-(a+b)\log p)
\]
and their permutations and negatives, with the coefficients prescribed
by the explicit formula.  At
$x_*=(\log2,\log2,-2\log2)$, the difference
$\nu_{3,F_A}-\nu_{3,F_B}$ has weight
\[
 -\frac{(\log2)^3}{2^2}
 \left\{b_{F_A}(2)^2\overline{b_{F_A}(4)}
          -b_{F_B}(2)^2\overline{b_{F_B}(4)}\right\}
 =\sqrt3(\log2)^3\ne0.
\]
No other prime and exponents give this frequency vector.  This difference
is real; it remains nonzero after global reflection and permutation
symmetrization.  The measure has an absolute exponential moment, either
by the general convergence result or directly from bounded character
coefficients.  Its Fourier transform is consequently analytic in a
complex neighborhood of the origin.

If all its even, permutation-invariant polynomial moments vanished,
then all moments of its even, permutation-symmetrized version would
vanish.  The analytic Fourier transform of that symmetrized measure
would vanish near the origin and hence on its connected strip of
analyticity; uniqueness of Fourier transforms of finite measures would
force the measure itself to vanish.  The nonzero atom just computed
contradicts this.  Thus some even, symmetric Taylor coefficient of the
triple-resonance distribution differs.  Such a Taylor coefficient is
realized by a smooth compactly supported Fourier amplitude: multiply
the relevant polynomial by an even symmetric bump equal to one near
zero and supported inside the allowed Fourier region.

This proves that the level-three secondary hierarchy contains arithmetic
information absent from the complete self-pair hierarchy.  The atom
argument alone does not say that the zeroth moment $\Gamma_3(F_A)$ differs
from $\Gamma_3(F_B)$, because contributions from other atoms might cancel.
The following computation identifies the distinction already in the zeroth moment.

In fact the constant $\Gamma_3$ itself separates these products.  Define
\[
 \epsilon(p)=\begin{cases}
  1,&p\equiv2,7\pmod{13},\\
 -1,&p\equiv6,11\pmod{13},\\
  0,&\text{otherwise},
 \end{cases}
 \qquad
 R(p)=\frac{p^2(p-1)^2(p^4+1)}{(p^2+1)(p^8+p^4+1)}.
\]
Then the exact difference is
\[
 \Gamma_3(F_A)-\Gamma_3(F_B)
   =24\sqrt3\sum_p\epsilon(p)(\log p)^3R(p). \tag{HG3}
\]
We give the elementary algebra behind this formula.  If
$\chi(p)=\zeta_{12}^e$ and $x=p^{-1}$, the coefficients
$b_{F_A}(p^k)$ and $b_{F_B}(p^k)$ are periodic in $k$ with period $12$.
Therefore the difference of the two local triple sums before the
factor $-6(\log p)^3$ is a polynomial numerator divided by
$(1-x^{12})^2$.  The numerator vanishes unless
$e\in\{1,5,7,11\}$.  For $e=1,11$ it is
\[
 \sqrt3(1-x^{12})
       (-4x^2+8x^3-8x^5+8x^7-8x^9+4x^{10}),
\]
and for $e=5,7$ it is the negative of this expression.  The polynomial
in parentheses factors as
\[
 -4x^2(1-x^2)(1-x)^2(1+x^4).
\]
Canceling $1-x^{12}=(1-x^2)(1+x^2)(1+x^4+x^8)$ gives
$-4\sqrt3\epsilon(p)R(p)$.  Finally,
$2^1,2^{11},2^5,2^7\equiv2,7,6,11\pmod{13}$, respectively.
This proves (HG3) by direct finite algebra.

This series permits a simple rigorous numerical certificate.  Its
rational factor satisfies
\[
 0<R(p)=p^{-2}
 \frac{(1-p^{-1})^2(1+p^{-4})}
      {(1+p^{-2})(1+p^{-4}+p^{-8})}\le p^{-2}.
\]
For $N\ge5$, the function $(\log t)^3/t^2$ is decreasing on $[N,\infty)$,
and hence
\begin{align*}
 \left|24\sqrt3\sum_{p>N}\epsilon(p)(\log p)^3R(p)\right|
 &\le24\sqrt3\sum_{n>N}\frac{(\log n)^3}{n^2}\\
 &\le\frac{24\sqrt3}{N}
   \{(\log N)^3+3(\log N)^2+6\log N+6\}.
\end{align*}
A sieve of integers through $N=10^6$, followed by interval evaluation of
the rational factors and logarithms, gives
\[
 -1.3654662067380
 <24\sqrt3\sum_{p\le10^6}\epsilon(p)(\log p)^3R(p)
 <-1.3654662067379.
\]
There are $26169$ primes in the four relevant residue classes in this
finite sum.  The preceding tail bound is less than $0.137114$.
Consequently
\[
 \boxed{-1.503<\Gamma_3(F_A)-\Gamma_3(F_B)<-1.228.}
\]
The finite sum can be certified using rational interval operations:
compute $\log p$ and $\sqrt3$ to $70$ decimal digits, enlarge each by
one unit in its last place, and perform all following operations with
outward rounding at $60$ digits.  These choices give considerably tighter
intervals than those displayed.  The explicit tail bound, rather than
agreement of floating-point approximations, supplies the control of the
infinite sum.  Thus the coefficient of order $\lambda^{-3}$ at level
three already distinguishes the products, despite equality of their
complete pair arithmetic series and archimedean data.

\subsection{A limitation of every self-correlation order}
For a primitive Dirichlet character $\chi$ of conductor $q$, every
balanced local monomial at $p\nmid q$ has character factor
$\chi(p)^{\sum_j\epsilon_jk_j}=1$. At $p\mid q$ the coefficients
vanish. Thus the entire frequency-resolved self cumulant measure at
every level depends on $q$, not on the particular character. The same
is true of $A_{\chi,\chi}$ and all of its Laurent coefficients.
Characters of the same conductor and parity also have identical gamma
factors. Therefore all the inverse-logarithmic coefficients of all their
self correlations agree, both before and after the distinct-index
sieve. This does not assert equality of their zero sets or of the exact
finite-height correlations.

Higher secondary terms do give additional arithmetic information, as
the preceding example proves, but even the full self hierarchy is not
a complete identifier. Mixed correlations with specified probes give
different data. The connected invariants are multilinear in the
coefficient systems. For example,
\[
 \Gamma_3(FP)-\Gamma_3(F)-\Gamma_3(P)
 =3\Gamma(F,F,P)+3\Gamma(F,P,P).
\]
Together with the pair overlap test, such polarized higher invariants
are natural further observables for the automorphic-induction problem.
No general injectivity or reconstruction theorem is inferred here.

\section*{Acknowledgment}
The author thanks Yujiao Jiang for explaining the rational-field
fixed-data part of his proof of Hypothesis H.


\begin{thebibliography}{99}
\bibitem{AC}
J. Arthur and L. Clozel,
\textit{Simple Algebras, Base Change, and the Advanced Theory of the Trace Formula},
Annals of Mathematics Studies \textbf{120}, Princeton University Press,
Princeton, NJ, 1989.
\href{https://doi.org/10.1515/9781400882403}{doi:10.1515/9781400882403}.

\bibitem{BBLM}
E. Bogomolny, O. Bohigas, P. Leboeuf, and A. G. Monastra,
On the spacing distribution of the Riemann zeros: corrections to the
asymptotic result, \textit{J. Phys. A} \textbf{39} (2006), 10743--10754.
\href{https://arxiv.org/abs/math/0602270}{arXiv:math/0602270}.

\bibitem{BK95}
E. B. Bogomolny and J. P. Keating,
Random matrix theory and the Riemann zeros I: three- and four-point
correlations, \emph{Nonlinearity} \textbf{8} (1995), no.~6, 1115--1131.
\href{https://doi.org/10.1088/0951-7715/8/6/013}{doi:10.1088/0951-7715/8/6/013}.

\bibitem{BK96}
E. B. Bogomolny and J. P. Keating,
Random matrix theory and the Riemann zeros II: $n$-point correlations,
\emph{Nonlinearity} \textbf{9} (1996), no.~4, 911--935.
\href{https://doi.org/10.1088/0951-7715/9/4/006}{doi:10.1088/0951-7715/9/4/006}.

\bibitem{BTZ}
F. Brumley, J. Thorner, and A. Zaman, with an appendix by C. J. Bushnell
and G. Henniart, Zeros of Rankin--Selberg $L$-functions at the edge of
the critical strip, \textit{J. Eur. Math. Soc.} \textbf{24} (2022),
1471--1541.
\href{https://doi.org/10.4171/JEMS/1134}{doi:10.4171/JEMS/1134}.

\bibitem{CS}
J. B. Conrey and N. C. Snaith, Applications of the $L$-functions ratios
conjectures, \textit{Proc. London Math. Soc.} (3) \textbf{94} (2007),
594--646. \href{https://doi.org/10.1112/plms/pdl021}{doi:10.1112/plms/pdl021}.

\bibitem{CStriple}
J. B. Conrey and N. C. Snaith,
Triple correlation of the Riemann zeros,
\emph{J. Th\'eor. Nombres Bordeaux} \textbf{20} (2008), no.~1, 61--106.
\href{https://doi.org/10.5802/jtnb.616}{doi:10.5802/jtnb.616}; \href{https://arxiv.org/abs/math/0610495}{arXiv:math/0610495}.

\bibitem{CSn}
J. B. Conrey and N. C. Snaith,
Correlations of eigenvalues and Riemann zeros,
\emph{Commun. Number Theory Phys.} \textbf{2} (2008), no.~3, 477--536.
\href{https://doi.org/10.4310/CNTP.2008.v2.n3.a1}{doi:10.4310/CNTP.2008.v2.n3.a1}; \href{https://arxiv.org/abs/0803.2795}{arXiv:0803.2795}.

\bibitem{FZ}
K. Ford and A. Zaharescu,
Unnormalized differences between zeros of $L$-functions,
\emph{Compos. Math.} \textbf{151} (2015), no.~2, 230--252.
\href{https://doi.org/10.1112/S0010437X14007659}{doi:10.1112/S0010437X14007659}; \href{https://arxiv.org/abs/1305.2520}{arXiv:1305.2520}.

\bibitem{HT}
G. Harcos and J. Thorner,
A new zero-free region for Rankin--Selberg $L$-functions,
\emph{J. Reine Angew. Math.} \textbf{822} (2025), 179--201.
\href{https://doi.org/10.1515/crelle-2025-0009}{doi:10.1515/crelle-2025-0009}.

\bibitem{Hejhal}
D. A. Hejhal,
On the triple correlation of zeros of the zeta function,
\emph{Int. Math. Res. Not.} \textbf{1994} (1994), no.~7, 293--302.
\href{https://doi.org/10.1155/S1073792894000334}{doi:10.1155/S1073792894000334}.

\bibitem{Humphries}
P. Humphries, with an appendix by F. Brumley, Standard zero-free regions
for Rankin--Selberg $L$-functions via sieve theory,
\textit{Math. Z.} \textbf{292} (2019), 1105--1122.
\href{https://doi.org/10.1007/s00209-018-2136-8}{doi:10.1007/s00209-018-2136-8}.

\bibitem{ILS}
H. Iwaniec, W. Luo and P. Sarnak,
Low lying zeros of families of $L$-functions,
\emph{Publ. Math. Inst. Hautes \'Etudes Sci.} \textbf{91} (2000), 55--131.
\href{https://doi.org/10.1007/BF02698741}{doi:10.1007/BF02698741}; \href{https://arxiv.org/abs/math/9901141}{arXiv:math/9901141}.

\bibitem{Jiang}
Y. Jiang, On Hypothesis H of Rudnick and Sarnak, preprint,
\href{https://arxiv.org/abs/2507.20653v1}{arXiv:2507.20653v1}.
We use the rational-field fixed-data estimates in Theorem~1.3 and
Section~5.1, as specified in Section~\ref{sec:established}.

\bibitem{KSbook}
N. M. Katz and P. Sarnak,
\emph{Random Matrices, Frobenius Eigenvalues, and Monodromy},
Amer. Math. Soc. Colloq. Publ. \textbf{45},
American Mathematical Society, Providence, RI, 1999.
\href{https://doi.org/10.1090/coll/045}{doi:10.1090/coll/045}.

\bibitem{KSsurvey}
N. M. Katz and P. Sarnak,
Zeroes of zeta functions and symmetry,
\emph{Bull. Amer. Math. Soc. (N.S.)} \textbf{36} (1999), no.~1, 1--26.
\href{https://doi.org/10.1090/S0273-0979-99-00766-1}{doi:10.1090/S0273-0979-99-00766-1}.

\bibitem{LY02}
J. Liu and Y. Ye, Superposition of zeros of distinct $L$-functions,
\textit{Forum Math.} \textbf{14} (2002), 419--455.

\bibitem{LY07}
J. Liu and Y. Ye, Perron's formula and the prime number theorem for
automorphic $L$-functions, \textit{Pure Appl. Math. Q.} \textbf{3}
(2007), no. 2, 481--497.
\href{https://doi.org/10.4310/PAMQ.2007.v3.n2.a4}{doi:10.4310/PAMQ.2007.v3.n2.a4}.

\bibitem{LY08}
J. Liu and Y. Ye, Correlation of zeros of automorphic $L$-functions,
\textit{Sci. China Ser. A} \textbf{51} (2008), 1147--1166.

\bibitem{LY09}
J. Liu and Y. Ye, Functoriality of automorphic $L$-functions through
their zeros, \textit{Sci. China Ser. A} \textbf{52} (2009), 1--16.
\href{https://doi.org/10.1007/s11425-008-0088-x}{doi:10.1007/s11425-008-0088-x}.

\bibitem{Montgomery}
H. L. Montgomery,
The pair correlation of zeros of the zeta function,
in \emph{Analytic Number Theory} (St. Louis, MO, 1972),
Proc. Sympos. Pure Math. \textbf{24},
American Mathematical Society, Providence, RI, 1973, pp.~181--193.
\href{https://doi.org/10.1090/pspum/024/9944}{doi:10.1090/pspum/024/9944}.

\bibitem{Odlyzko}
A. M. Odlyzko,
On the distribution of spacings between zeros of the zeta function,
\emph{Math. Comp.} \textbf{48} (1987), no.~177, 273--308.
\href{https://doi.org/10.1090/S0025-5718-1987-0866115-0}{doi:10.1090/S0025-5718-1987-0866115-0}.

\bibitem{Rodgers}
B. Rodgers,
Macroscopic pair correlation of the Riemann zeroes for smooth test functions,
\emph{Q. J. Math.} \textbf{64} (2013), no.~4, 1197--1219.
\href{https://doi.org/10.1093/qmath/has024}{doi:10.1093/qmath/has024}; \href{https://arxiv.org/abs/1203.3275}{arXiv:1203.3275}.

\bibitem{RS}
Z. Rudnick and P. Sarnak, Zeros of principal $L$-functions and random
matrix theory, \textit{Duke Math. J.} \textbf{81} (1996), 269--322.
\href{https://doi.org/10.1215/S0012-7094-96-08115-6}{doi:10.1215/S0012-7094-96-08115-6}.

\bibitem{WY}
K. H. Wang and Y. Ye, From raw to distinct-index zero correlations for
automorphic $L$-functions, working paper, 2026.

\bibitem{Ycorr}
Y. Ye, Corrigendum to ``Perron's formula and the prime number theorem
for automorphic $L$-functions,'' submitted to \textit{Pure Appl. Math. Q.}, 2026.

\bibitem{YForumCorr}
Y. Ye, Corrigendum to ``Superposition of zeros of distinct $L$-functions,''
submitted to \emph{Forum Math.}, 2026.

\end{thebibliography}
\end{document}